\documentclass[english,a4paper,11pt]{article} 
\usepackage{amsmath,amssymb}
\usepackage{graphicx,url}
\usepackage{babel}

\newcommand{\RR}{\mathbb{R}}

\newcommand{\ZZ}{\mathbb{Z}}
\newcommand{\CC}{\mathbb{C}}
\newtheorem{Theorem}{Theorem}

\begin{document}

\title{Elementary fractal geometry. \\  6. The dynamical interior of self-similar sets}
\author{Christoph Bandt\\ 
Institute of Mathematics, University of Greifswald,  
\url{bandt@uni-greifswald.de}}
\maketitle

\begin{abstract} On the one hand, the dynamical interior of a self-similar set with open set condition is the complement of the dynamical boundary. On the other hand, the dynamical interior is the recurrent set of the magnification flow.  For a finite type self-similar set, both boundary and interior are described by finite automata.  The neighbor graph defines the boundary. The neighborhood graph, based on work by Thurston, Lalley, Ngai and Wang, defines the interior. If local views are considered up to similarity, the interior obtains a discrete manifold structure, and the magnification flow is discretized by a Markov chain. This leads to new methods for the visualization and description of finite type attractors. 
\end{abstract}

\section{Introduction}  \label{intro}
{\bf Self-similar sets.}
A compact non-empty set $A$ in Euclidean $\RR^d$ is called self-similar if it is a finite union of subsets $A_j$ which are geometrically similar to $A,$ cf.~\cite{Mo}.  Using contractive similitudes $f_j$ and writing $A_j=f_j(A),$ Hutchinson has characterized these sets as compact solutions of the equation
\begin{equation} A=\bigcup_{j=1}^m f_j(A) \ .  \tag{1}
\label{hut}\end{equation}
The set $F=\{ f_1,...,f_m\}$ is called an iterated function system or IFS, and $A$ is the attractor of $F.$ A contractive similitude from Euclidean $\RR^d$ to itself fulfils $|f(x)-f(y)|=r_f|x-y|$ where the constant  $r_f<1$ is called the factor of $f.$ Unless explicitly stated, we assume that all maps in $F$ have the same factor $r.$ See \cite{BSS,Bar,BP,Fal,Fra} for basic concepts.\vspace{1ex}

{\bf The problem.}
It is not true, however,  that small subsets of a self-similar set $A$ are geometrically similar to $A.$ This holds only for the sets $f_w(A)=A_w$  with $w=w_1...w_n\in \{ 1,...,m\}^k$
and $f_w=f_{w_1}\cdots f_{w_k}.$  Such sets are very special and can hardly be recognized in magnifications, as shown in Figures \ref{figA}, \ref{figB}, \ref{figC} below.  These examples belong to the simplest self-similar sets.  They are modifications of the well-known Sierpi\'{n}ski triangle and fulfil the OSC. Nevertheless, their magnifications are quite diverse and do not look similar to each other or to $A.$ The purpose of the present note is to explain where the self-similarity can be found.\vspace{1ex}

{\bf OSC.}
An IFS $F$ fulfils the open set condition (OSC) if there is an open set $U$ such that the $f_j(U)$ are pairwise disjoint subsets of $U.$ The idea is that overlaps $f_i(A)\cap f_j(A)$ in this case must be small since they are contained both in the boundary of $f_i(U)$ and in the boundary of $f_j(U).$ Moran \cite{Mo} introduced the OSC in 1946 to prove that the uniform measure on $A$ is the normalized $\alpha$-dimensional Hausdorff measure. Thus from a measure-theoretic point of view, such attractors have a very homogeneous metric structure.

The OSC is crucial for the structure of $A$ \cite{BSS}. As Schief \cite{Sch} has noted, attractors without OSC are not truly self-similar since for pieces $A_w$ with level $k\to\infty ,$ the number of neighboring pieces of the same size will tend to infinity. However, one must not imagine that $U$ is a circle or triangle, a convex or just a connected set. The typical set $U$ will have infinitely many components, as pointed out in \cite{EFG1,Tet24}. This is why we can hardly identify the $f_w(A),$ even in the colored view of Figure \ref{figA}. 

\begin{figure}[h!t]
\begin{center}
\includegraphics[width=0.4\textwidth]{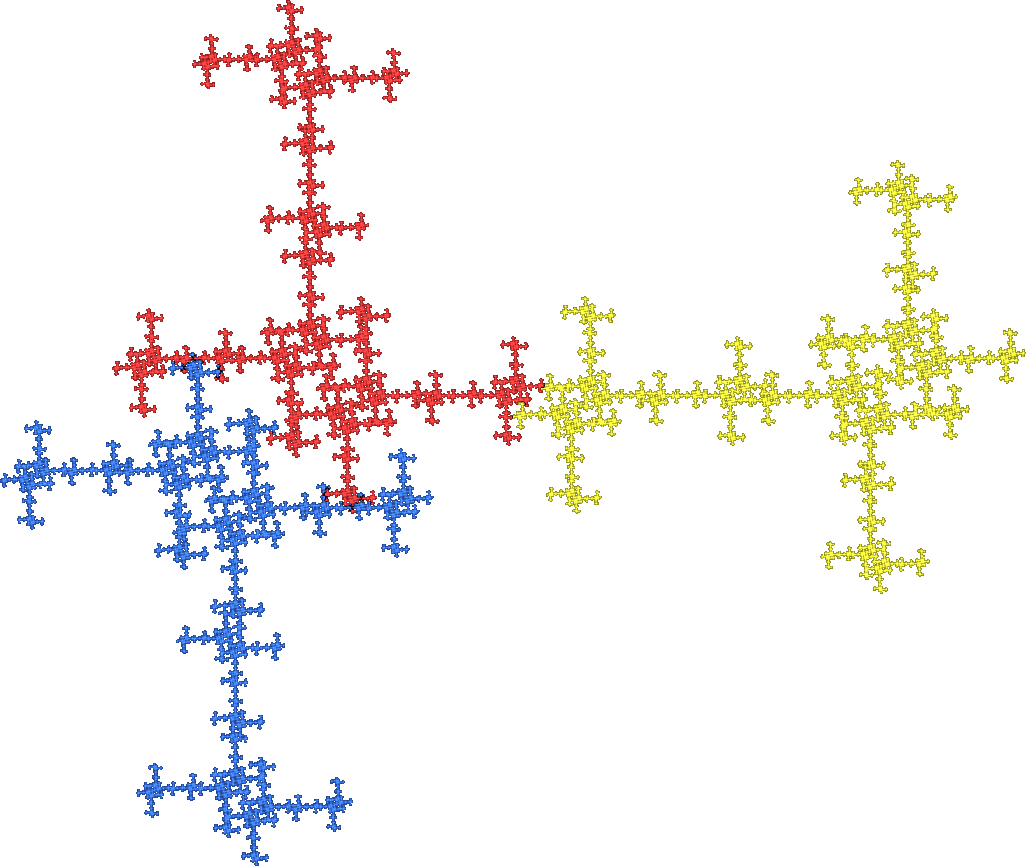} \qquad 
\includegraphics[width=0.4\textwidth]{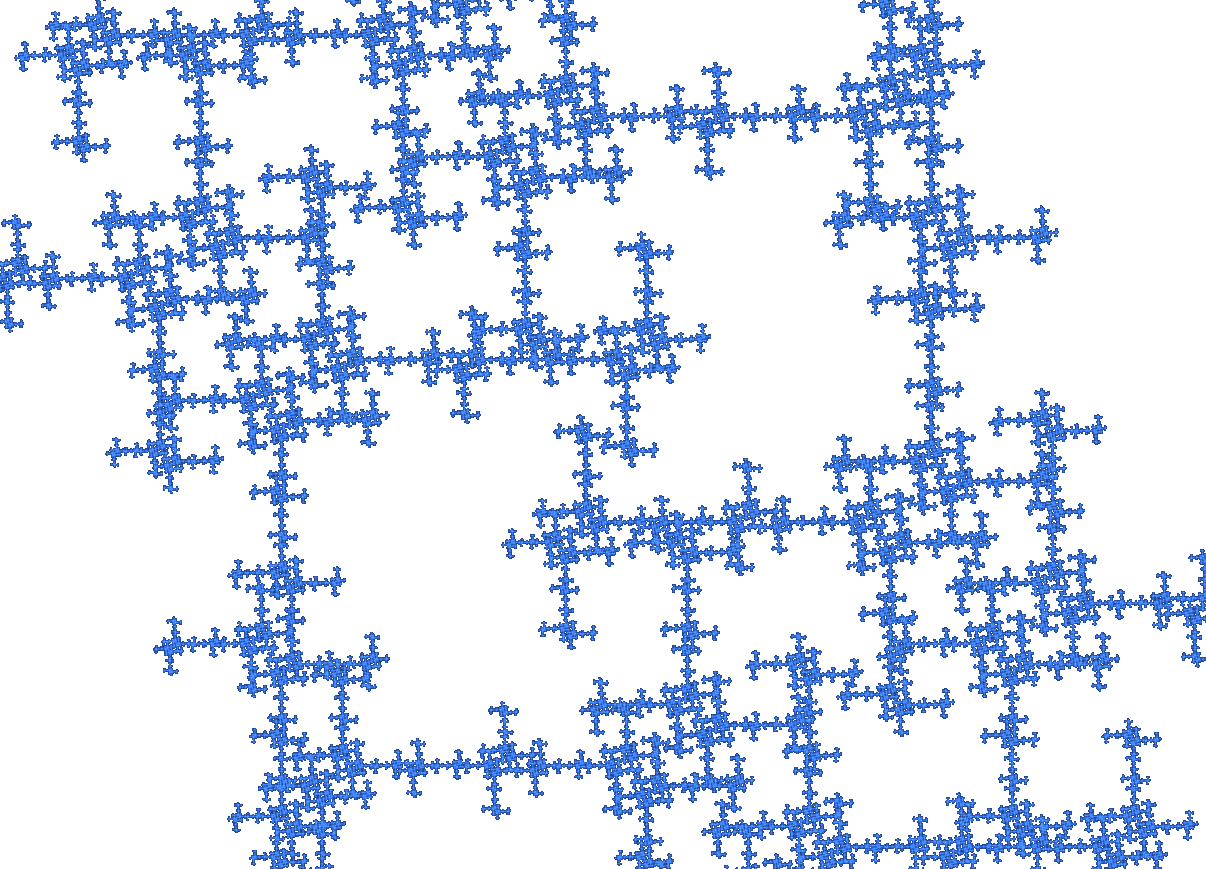}\vspace{1ex}\\
\includegraphics[width=0.4\textwidth]{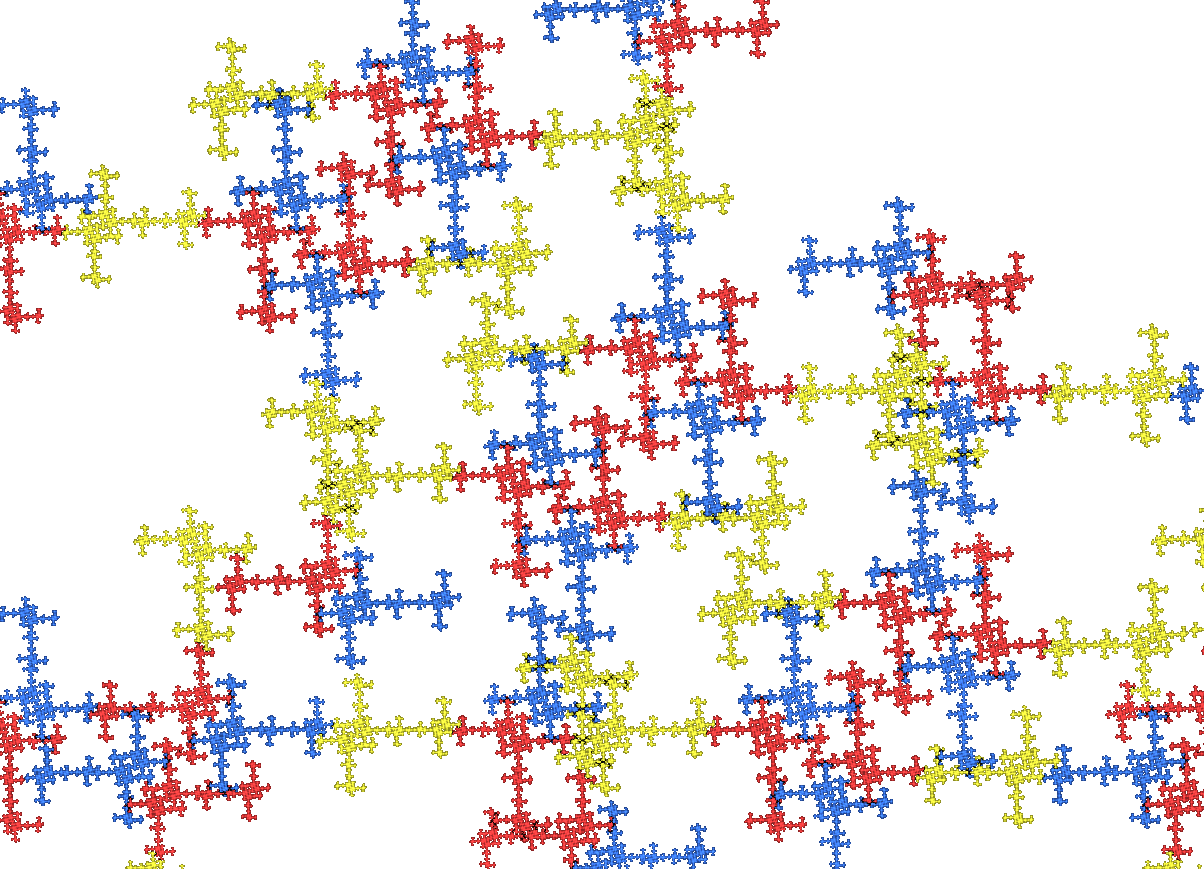} \qquad 
\includegraphics[width=0.4\textwidth]{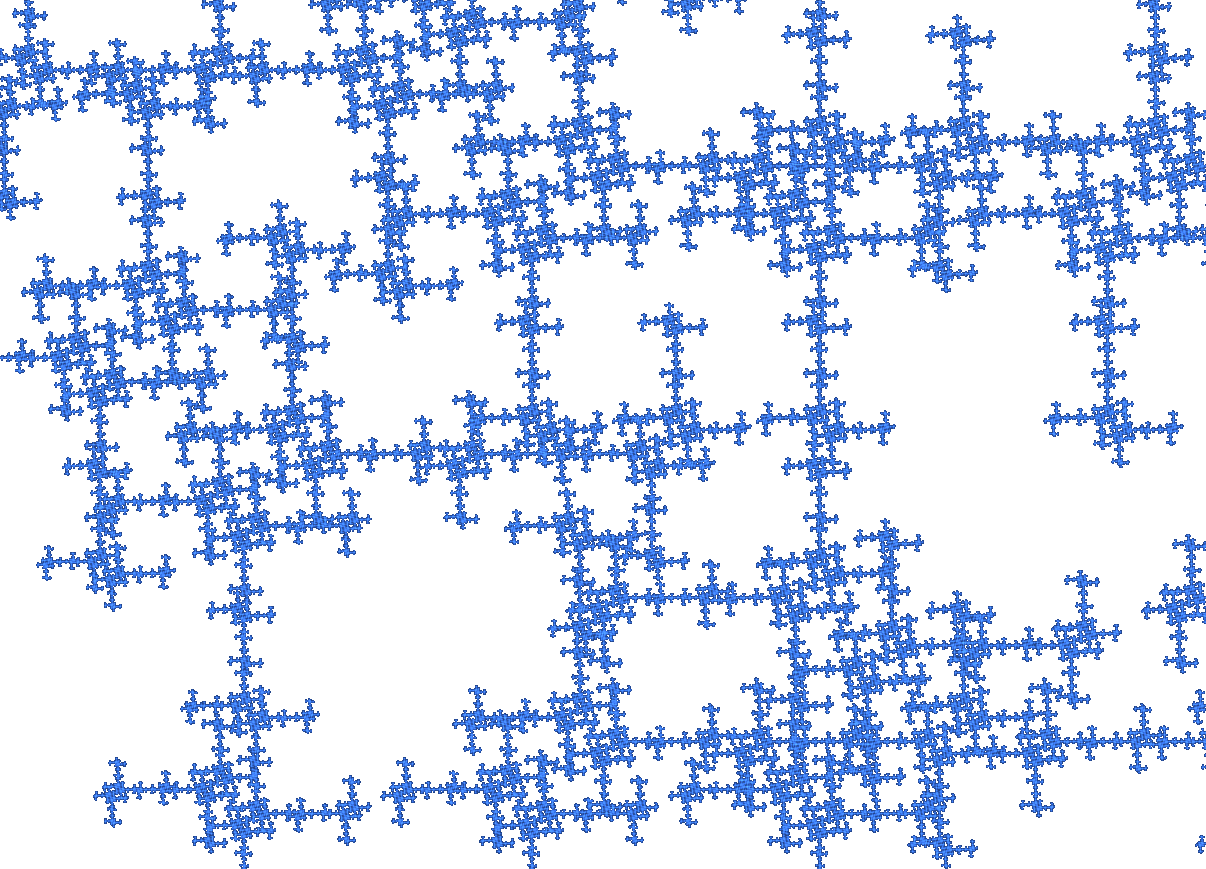}
\end{center}
\caption{The self-similar set A with three magnifications. In the second one, pieces are colored with respect to level 12. The IFS data are in Table \ref{tab1}.}
\label{figA}      
\end{figure}

{\bf Finite type.} When we strengthen the OSC, attractors become tractable. The IFS $\{ f_1,...,f_m\}$ is of finite type if there are only finitely many mappings $h=f_v^{-1}f_w$ with $v,w\in\{ 1,...,m\}^k$ for arbitrary $k$ such that $A\cap h(A)\not=\emptyset ,$ that is, $f_v(A)\cap f_w(A) \not=\emptyset .$ In other words, there are only finitely many configurations of two intersecting pieces of the same size, up to similarity. This implies that up to similarity, the configurations of two copies $f_v(U), f_w(U)$ with a common accumulation point are also finite in number. Thus both the attractor $A$ and also the open set $U$ with its infinite number of components must have a very special structure - for every admissible choice of $U.$ {In this paper, we assume both OSC and finite type.} That is, there are no exact overlaps $A_v=A_w$ \cite{EFG3}.\vspace{1ex}

{\bf Automata.} The advantage of the finite type property is that all topological properties of the attractor $A$ can be formally described by a finite automaton, called the neighbor graph \cite{EFG1,EFG3,EFG5}. Thus the geometrical structure of $A$ becomes accessible to computer work, without any reference to eyesight or geometric vision. Actually we need the computer to study even such simple examples like Figure \ref{figA} where the automaton has 23 states. 

The purpose of this paper is to derive from the neighbor graph another automaton, the neighborhood graph.  This concept is connected with ideas of Thurston \cite[Section 9]{T89}, Lalley \cite{Lal97}, and Ngai and Wang \cite{NW}.  In one-dimensional overlapping constructions it was used in various papers of Feng, K.E.~Hare, K.G.~Hare, Rutar, and others  \cite{EFG5,Feng16,HR22,rutar23}. Its general relevance has never been pointed out, however.  We show that the neighborhood graph of a finite type attractor $A$ is the appropriate tool to describe all magnifications of $A,$ and to navigate through their structure.
\vspace{1ex}

\begin{figure}[h!b]
\begin{center}
\includegraphics[width=0.4\textwidth]{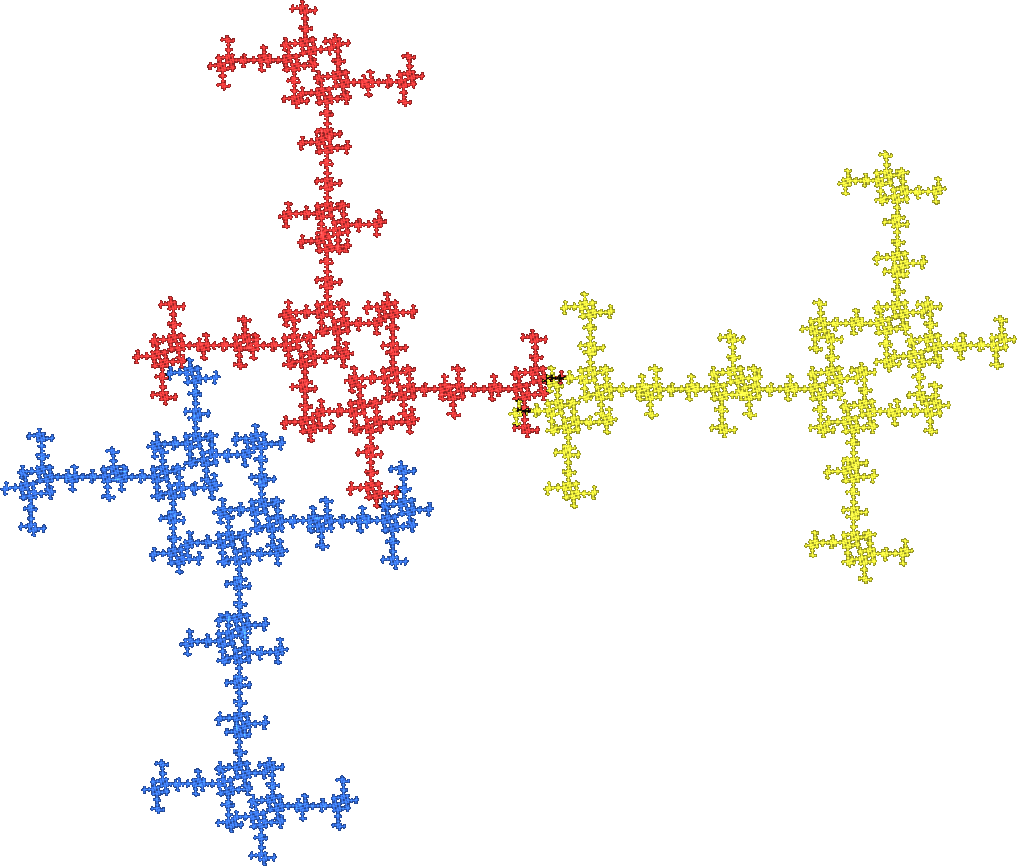} \qquad 
\includegraphics[width=0.4\textwidth]{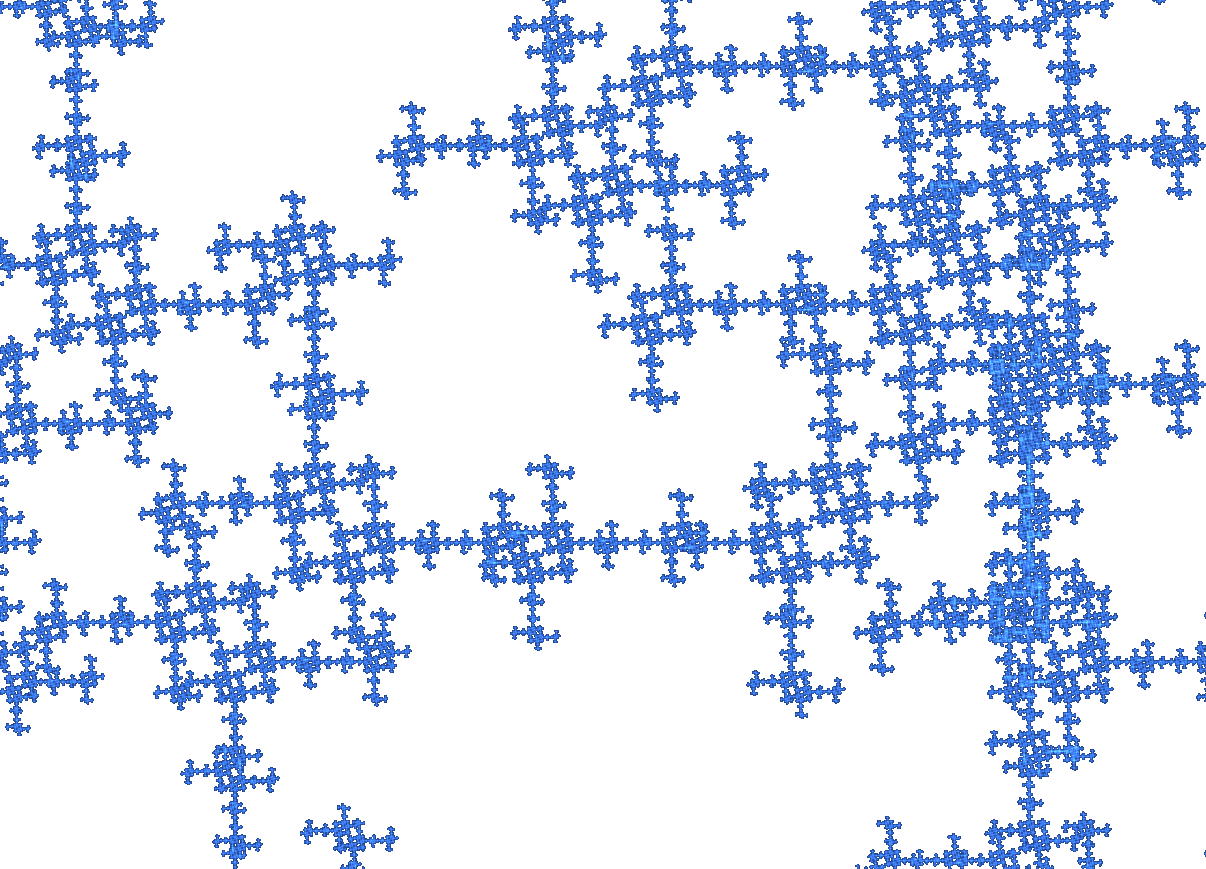}\vspace{1ex}\\
\includegraphics[width=0.4\textwidth]{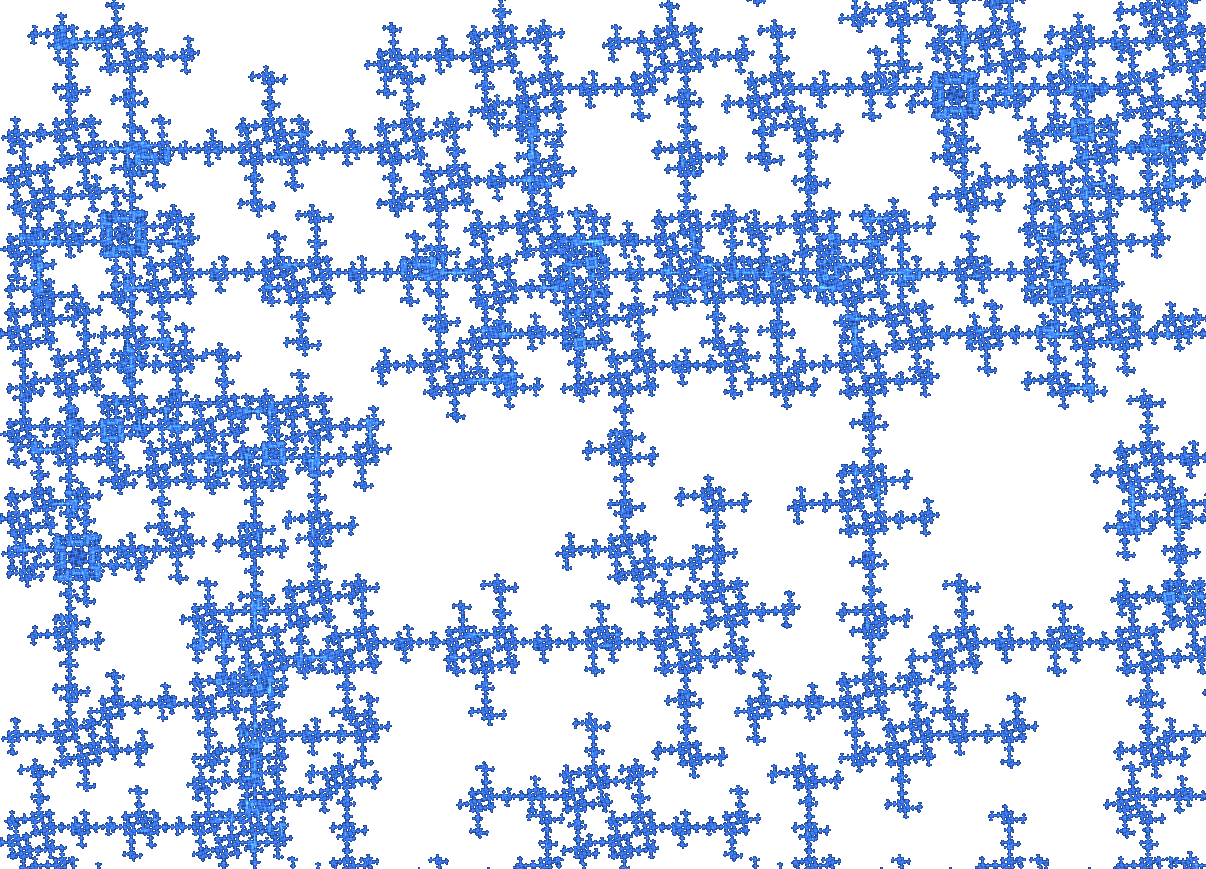} \qquad 
\includegraphics[width=0.4\textwidth]{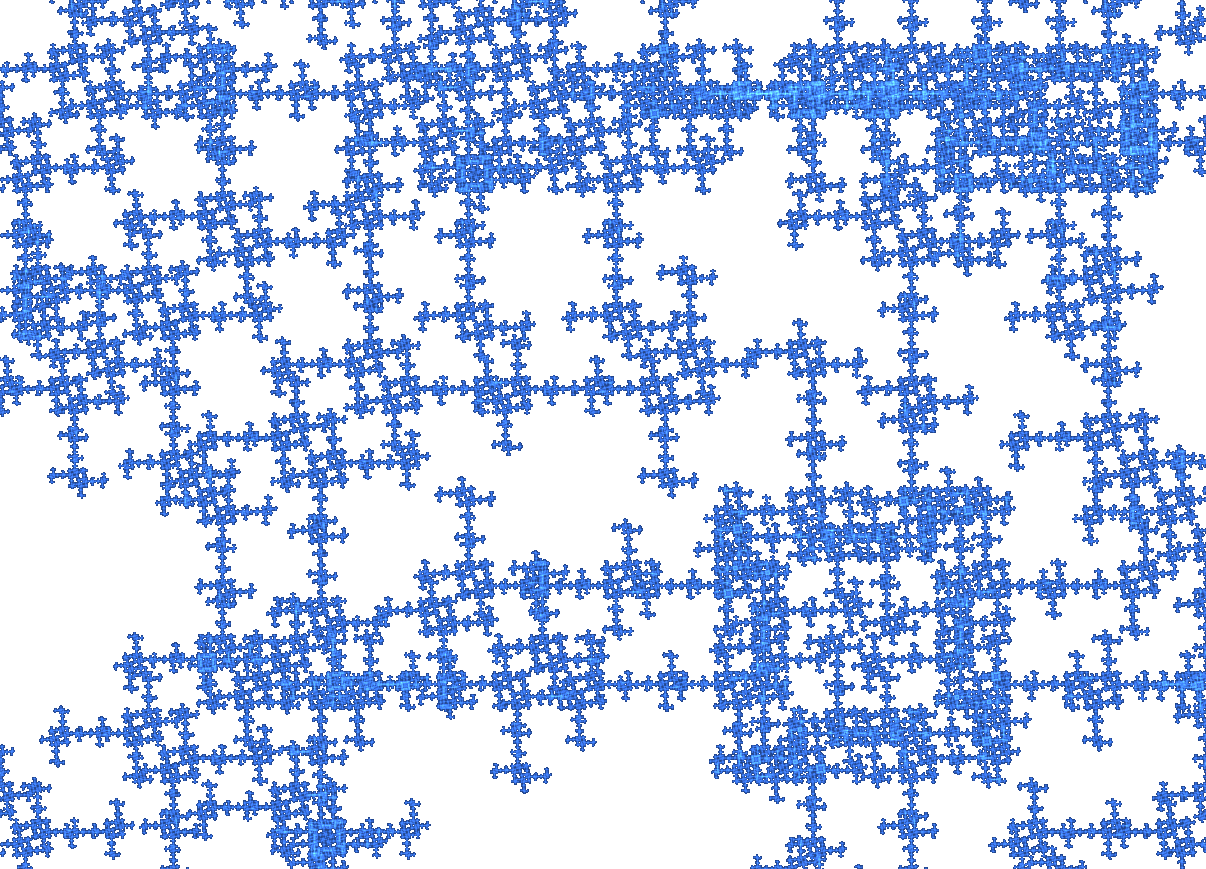} \vspace{1ex}\\
\includegraphics[width=0.4\textwidth]{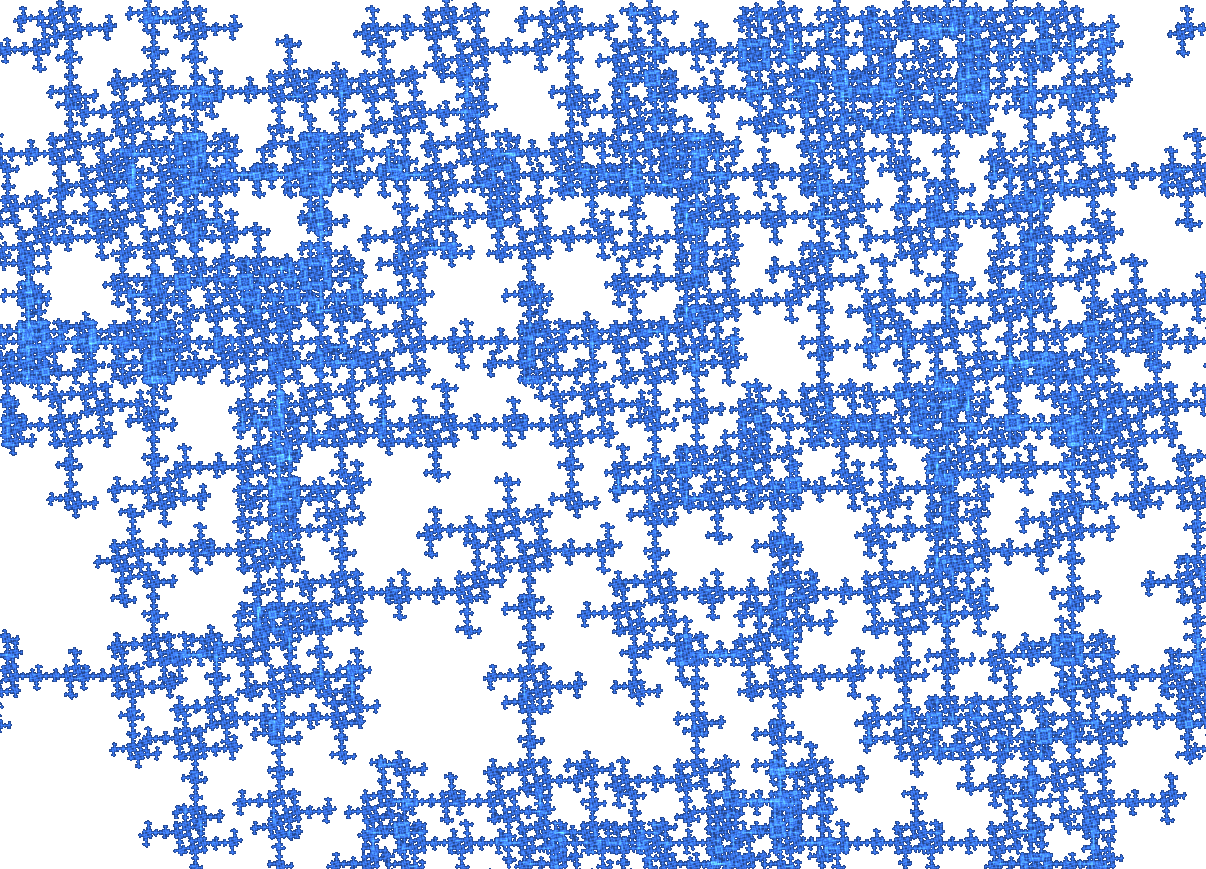} \qquad 
\includegraphics[width=0.4\textwidth]{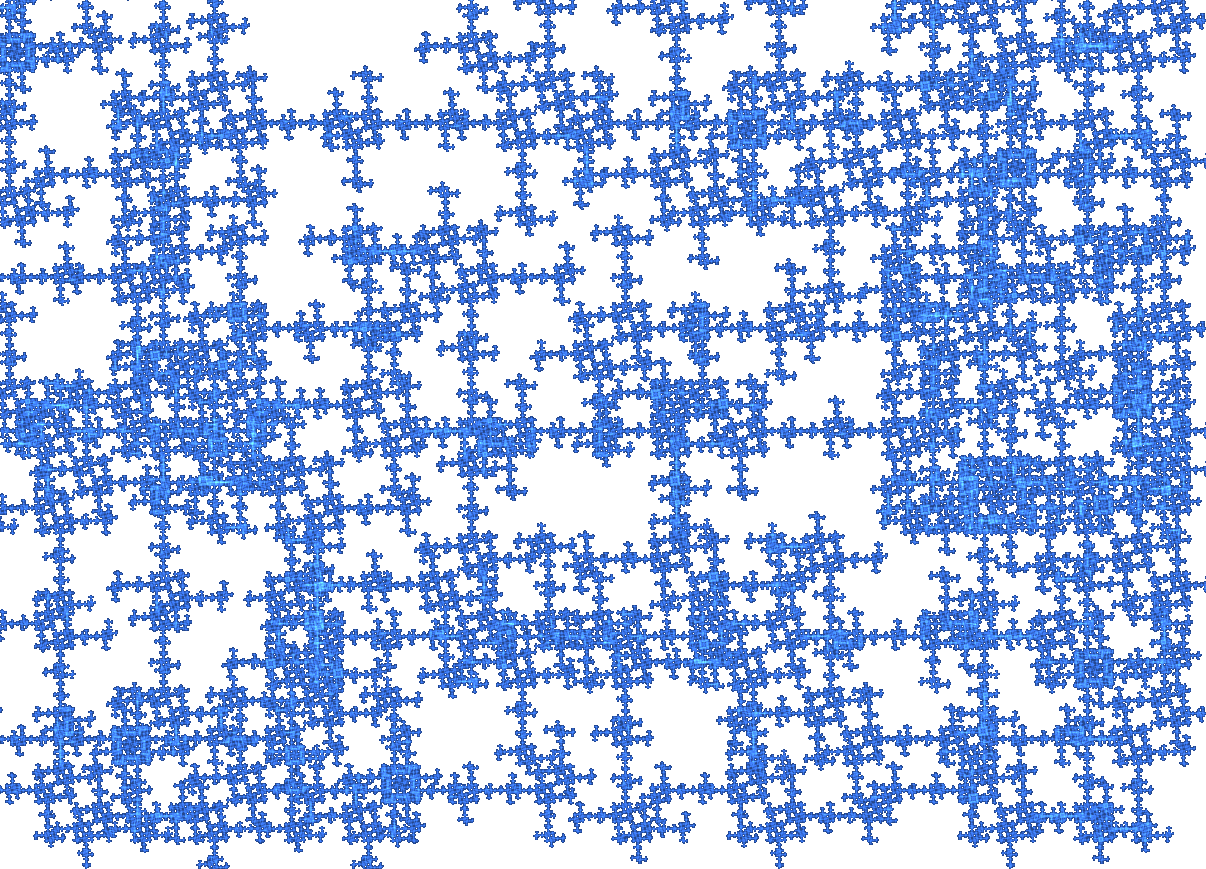}
\end{center}
\caption{The self-similar set B with a greater variety of magnifications. }
\label{figB}      
\end{figure}

{\bf Contents of the paper.} We start with some motivating examples. In Section \ref{conc} we formally introduce the neighbor graph and the boundary of a self-similar set.
The main and completely new part of the paper is Section \ref{inte}, where we state the definitions of dynamical interior and neighborhood graph.  We present a simple algorithm which derives the neighborhood graph from the neighbor graph. In Section \ref{out} we establish relations with the magnification flow and introduce the two-sided Markov chain associated with the neighborhood graph.  We discuss new methods for visualization and description of attractors, as well as problems for further research.

\section{Three examples}  \label{exa}
Using Mekhontsev's package IFStile \cite{M}, we study some IFS with three mappings and contraction factor $\frac12$ in the complex plane $\CC .$   The Sierpi\'{n}ski triangle with vertices $2v_1, 2v_2,$ and $2v_3$ is given by $f_j(z)=\frac{z}{2} +v_j, \ j=1,2,3.$ To get another structure, we have to include rotations in the IFS mappings. A $180^o$ rotation, obtained by replacing $\frac{z}{2}$ with $\frac{-z}{2}$ in some of the maps, will lead only to three other well-known examples.  The next trial is a $90^o$ rotation, represented by multiplication with $\pm i.$  This already leads to a huge family of OSC examples, studied in \cite{EFG1}. All attractors have Hausdorff dimension $\frac{\log 3}{\log 2}\approx 1.58,$ and they look a bit squarish.  We decided to take the following subfamily of IFS. 
\begin{equation}
f_1(z)=\frac{iz}{2} +v_1, \ f_2(z)=\frac{-z}{2} +v_2, \  f_3(z)=\frac{-iz}{2} +v_3\, . 
\tag{2} \label{IFS3}\end{equation}
That is, we downloaded IFStile, wrote the text file template.aifs below, clicked the file to open in IFStile, and clicked the binocular button to start the search. Mekhontsev's program is matrix-based, so $i$ had to be replaced by the matrix $s={0\, -1\choose 1\ \ 0}$ in dimension 2. It randomly generates the $v_j,$ tries to determine the neighbor graph for each example, and keeps only examples which are finite type and OSC. After three minutes you have about 100 attractors similar to those in our figures. You can order them in the list and watch their magnifications.

\verb|@F |

\verb|$dim=2 |

\verb|s=[0,-1,1,0] |

\verb|&T=$vector() |

\verb|f0=T*2^-1*s |

\verb|f1=T*2^-1*-1 |

\verb|f2=T*2^-1*s^3 |\vspace{-2ex}

\begin{verbatim}   A=(f0|f1|f2)*A\end{verbatim}

\begin{figure}[h!t]
\begin{center}
\includegraphics[width=0.4\textwidth]{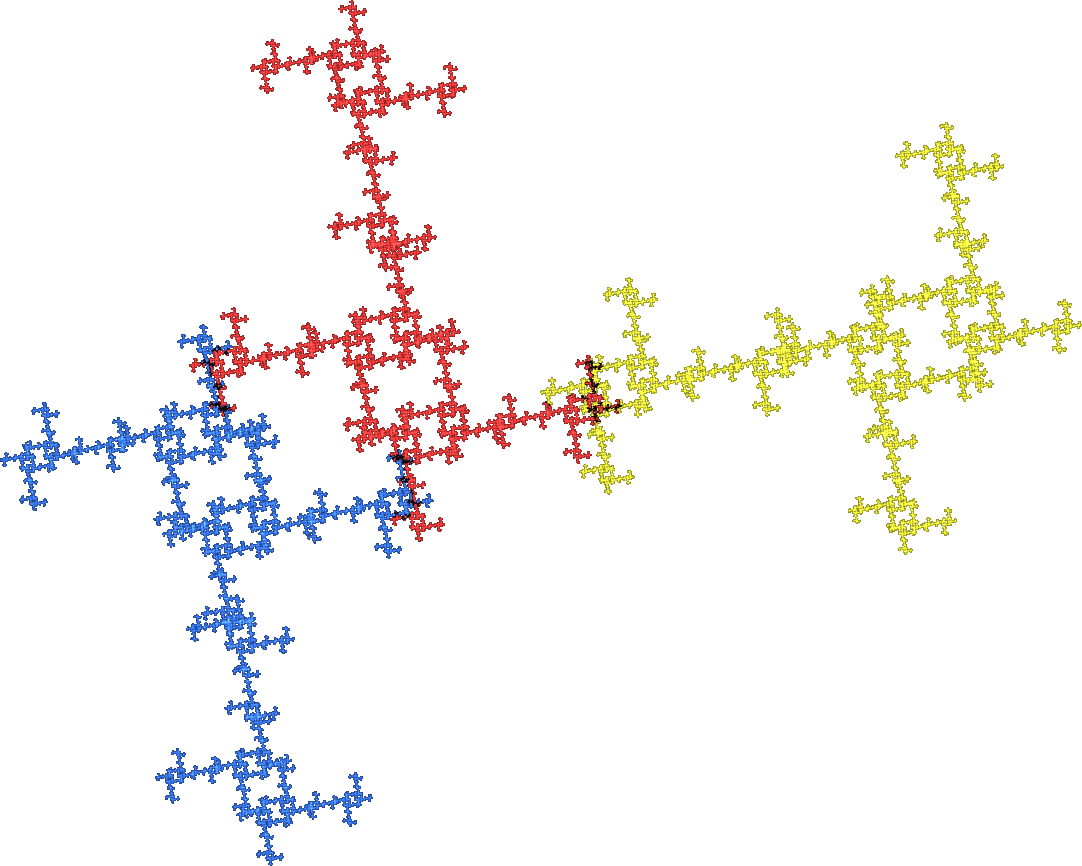} \qquad 
\includegraphics[width=0.4\textwidth]{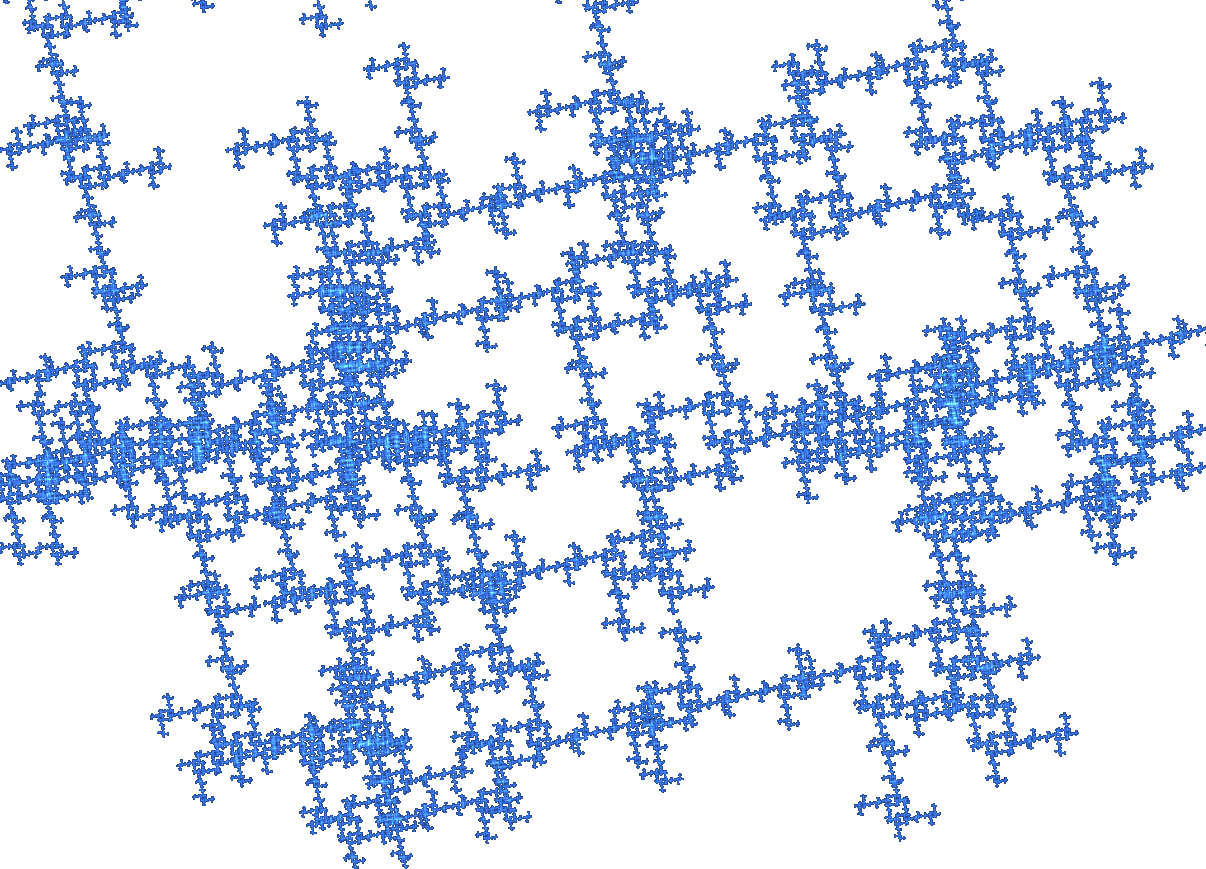}\vspace{1ex}\\
\includegraphics[width=0.4\textwidth]{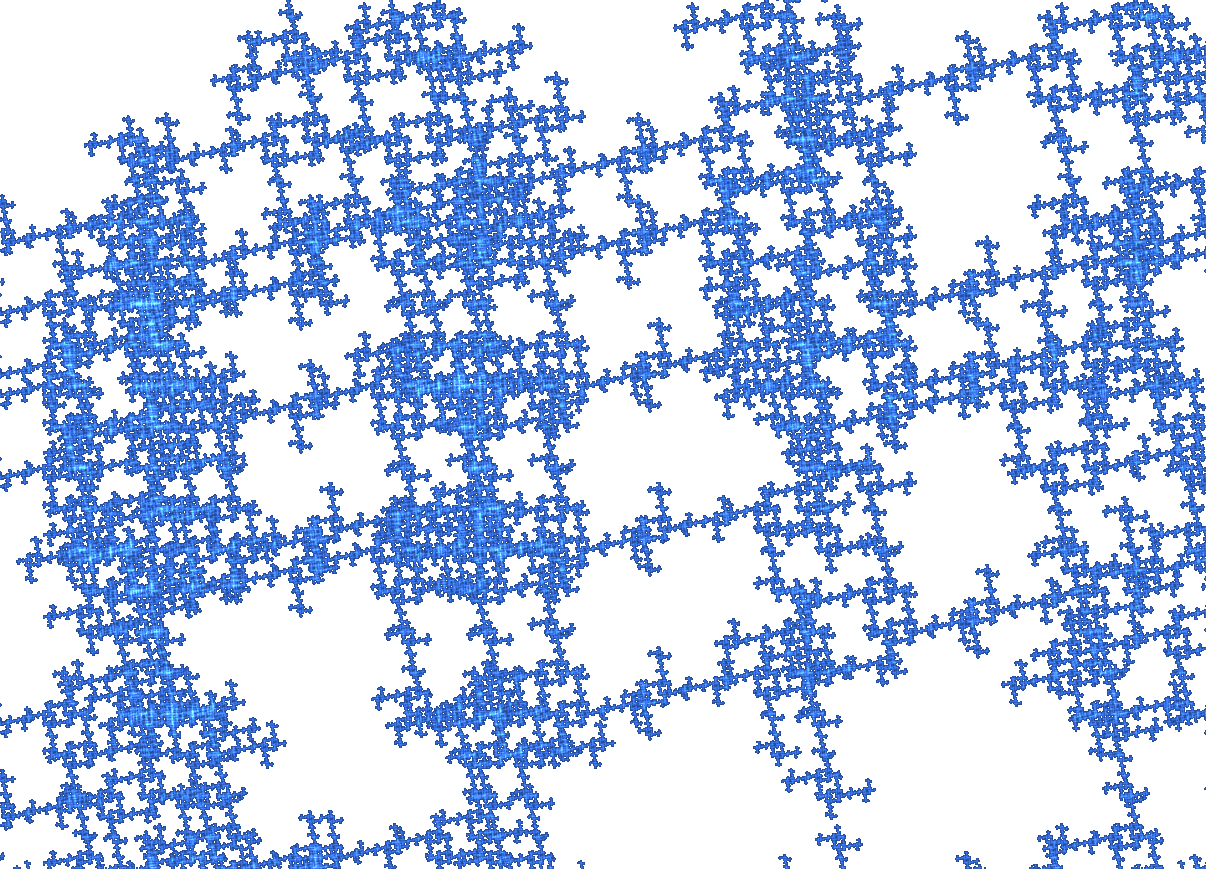} \qquad 
\includegraphics[width=0.4\textwidth]{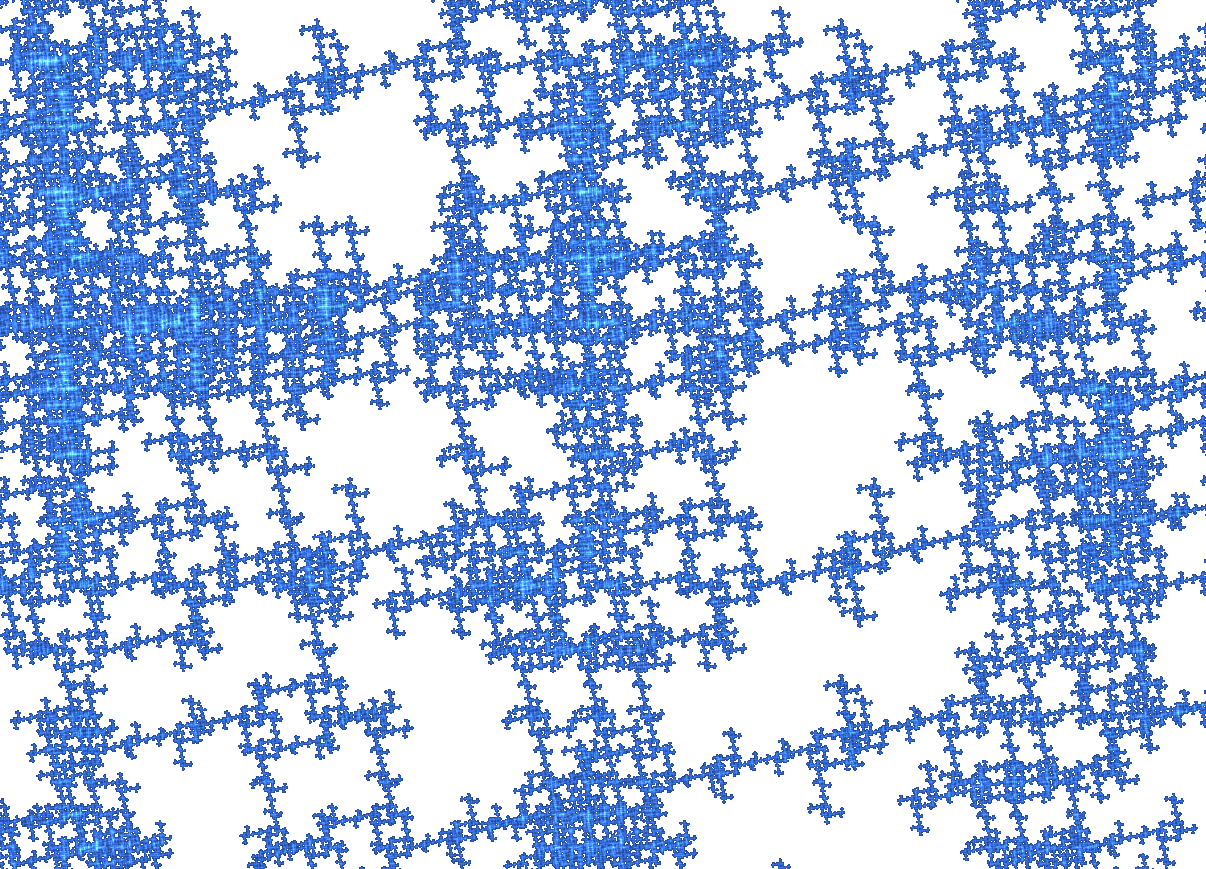}\vspace{1ex}\\
\includegraphics[width=0.4\textwidth]{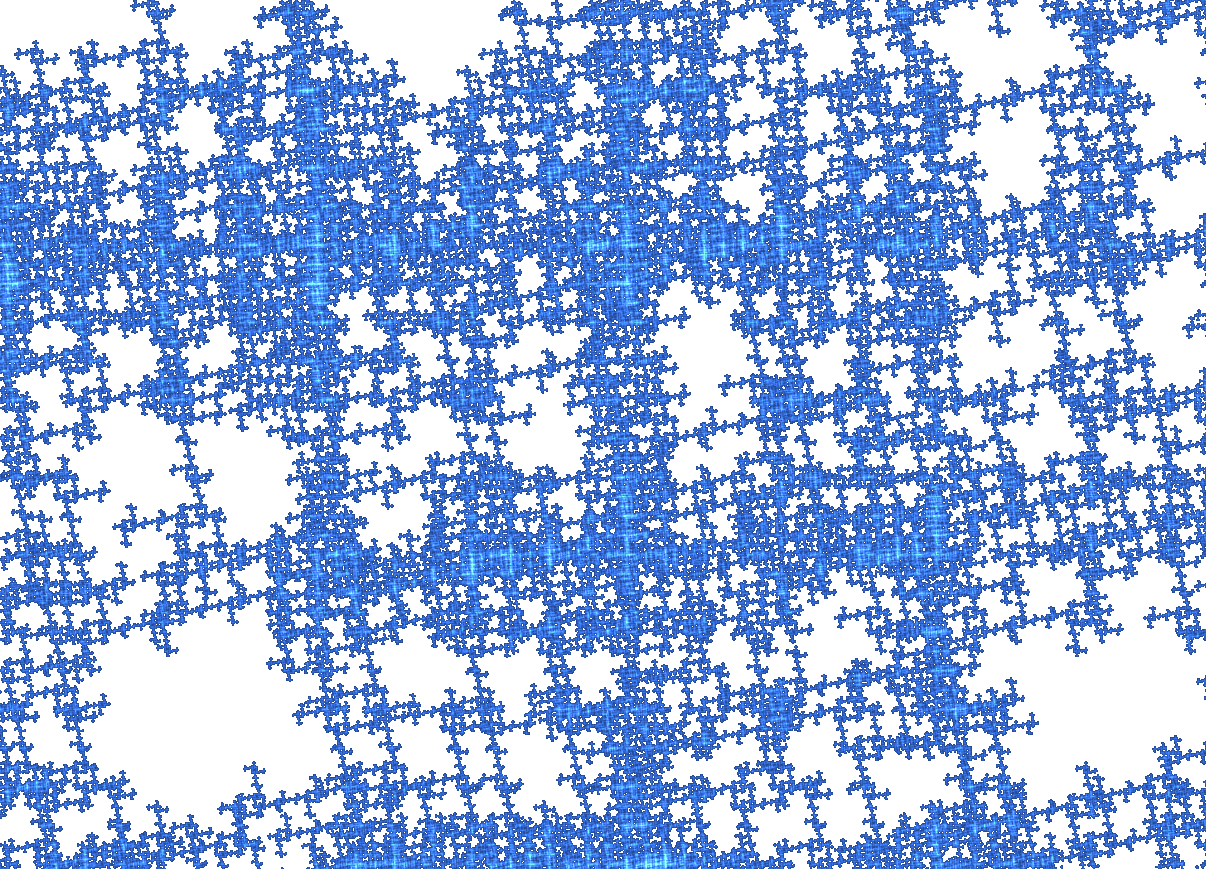} \qquad 
\includegraphics[width=0.4\textwidth]{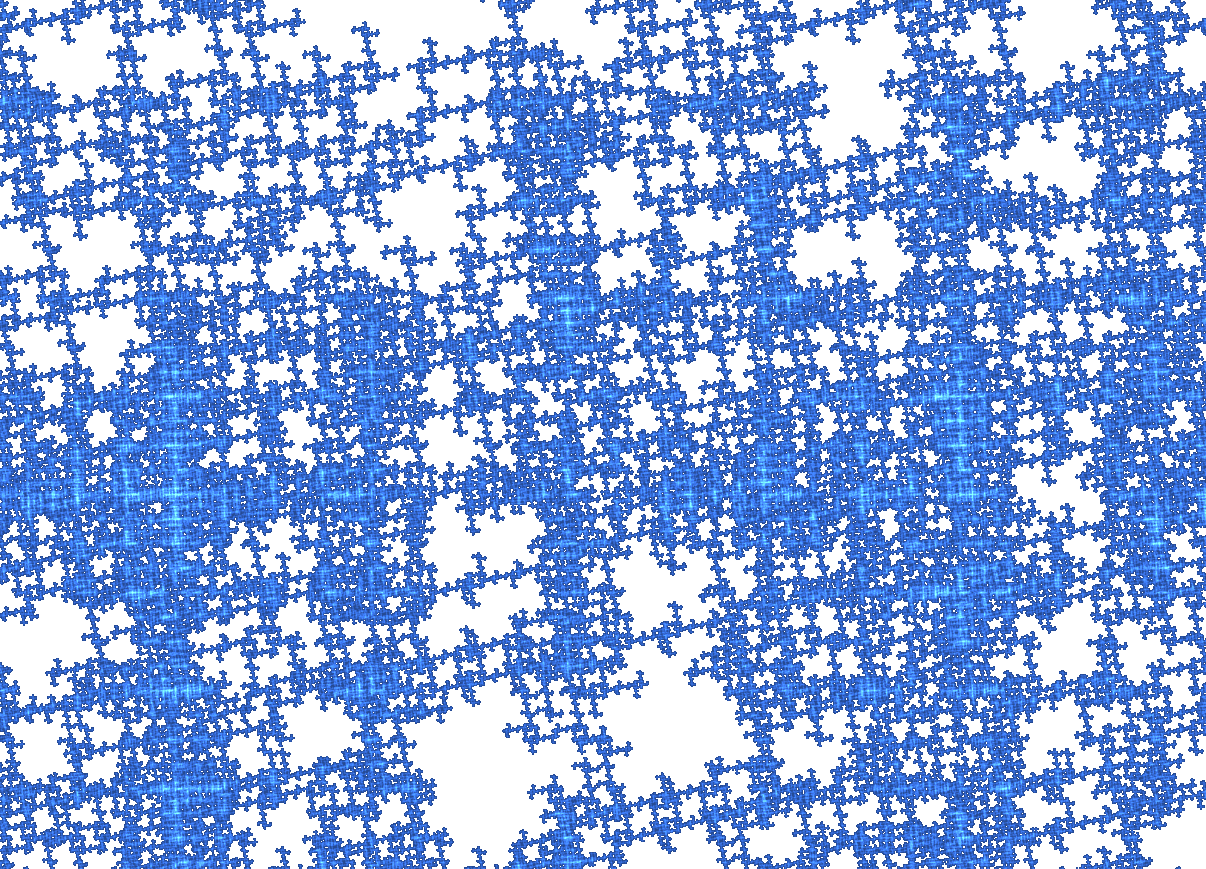}
\end{center}
\caption{The self-similar set C contains a lot of structural details. }
\label{figC}      
\end{figure}

The choice of three different rotations was motivated by the wish to have no line segments in the attractors. Actually, this is the case in our figures since no magnification contained lines which lead from one side of the window to the opposite one. We also chose the attractors to be connected, and to look fairly similar in the global view. 

Our goal was to show how much compact sets like $A,B,C$ differ in their local views. Investigating successive magnifications, we always chose some of the most dense or complex parts for the next view.  In Figure \ref{figA}, a few magnifications were sufficient to get an idea of the structure, even though each local view looked different. In Figure \ref{figB}, there are very dense parts, and the five views do not yet explain all the structural details. And in Figure \ref{figC} some of the views are still denser, and it would take a long time to explore all details.

This experiment implies that it is not enough to consider attractors of IFS just as compact sets.  They have an interior structure.  We must study local views.  It is like a vacation to an island. First you have the global map. Later you take lots of photos which provide much more information.  There are variations of the interior structure - in complexity, in density, in homogeneity - which are not expressed by Hausdorff dimension. Even the simplest self-similar sets have features which have never been studied. Are they truly self-similar? These questions will now be addressed more rigorously.

\begin{table}[h!t] \label{tab1}  
\begin{center}     
\begin{tabular}{|l|c|l|c|c|c|}
\hline
Attractor & $v_1, v_2, v_3$ & Boundary Dim & Neighbors & Neighborhoods\\
\hline
A & $i,\ 2+i,\ -1$ & 0.76 & 23 & 88\\
B & $2i,\ 3+2i,\ -2$ & 1.16 & 82 & 6291\\
C & $i,\ 5+i,\ -2-i$ & 1.21 & 186 & 168974\\
\hline
\end{tabular}
\end{center}
\caption{Data of the three examples. The IFS act on $\CC$ as $f_1(z)=\frac{iz}{2} +v_1, \ f_2(z)=\frac{-z}{2} +v_2,$ and $f_3(z)=\frac{-iz}{2} +v_3. $ All attractors have Hausdorff dimension $\frac{\log 3}{\log 2}\approx 1.58.$}
\end{table}

\section{Dynamical boundary and neighbor graph}  \label{conc}
{\bf Concept of boundary.}
The boundary of a set $A$ in a topological space $X$ consists of all points which can be approached from both $A$ and its complement: ${\rm bd}\, A= {\rm cl}\, A \cap {\rm cl}\, (X\setminus A)\ $ where cl denotes topological closure.  The boundary is the complement of the interior of $A.$ In particular, if $A$ has no interior points, which is the typical case for self-similar sets in $\RR^d,$ then $A$ coincides with its boundary. 

The dynamical boundary of a self-similar set $A$ was defined by M.~Moran \cite{Mo1}. It refers to the possibility of extending the self-similar construction of $A$ outwards, treating $A$ in the same way as a small piece $A_w$ together with all its intersecting pieces $A_v$.  The sets 
\begin{equation}  A\cap f^{-1}_w(A_w\cap A_v)= A\cap f^{-1}_wf_v(A)\ \mbox{ with } v,w\in \{ 1,...,m\}^k, \ k=1,2,...    \tag{3} \label{bdsets}
\end{equation}
are the boundary sets of $A.$ This is where $A$ touches a neighboring piece in a possible outward construction.  The closure of the union of all boundary sets is the dynamical boundary $B(A)$ of $A.$ 
This definition and the following theorem hold for arbitrary self-similar sets, without assuming equal contraction factors, OSC, or finite type. Part a) was shown in \cite[Lemma 2.1, (iii)]{Mo1} in a different way. Attractors with non-empty interior in c) where discussed in greater generality in \cite[Proposition 2.8]{Mo1}. It is still not known whether the OSC must hold if $A\not= B(A)$ (the proof of \cite[Theorem 2.3]{Mo1} contained a gap). 

\begin{Theorem}{\rm\cite{Mo1,BHR}} \label{bound}\\
Let $A$ be the self-similar set in $\RR^d$ with respect to similitudes $f_1,...,f_m.$
\begin{enumerate}
\item[a)]  The dynamical boundary $B(A)$ is disjoint to every open set $U$ which satisfies the OSC.  If $A=B(A),$ the OSC is not true.
\item[b)]  If the OSC holds, there exists an open set $U$ such that  $B(A)=A\setminus U.$
\item[c)]  Assume the OSC. Then $B(A)$ agrees with the topological boundary if and only if $A$ has non-empty interior.
\end{enumerate}
\end{Theorem}

{\it Proof. } 
a) The OSC says that $f_w(U)\cap f_v(U)=\emptyset$ whenever $w_1\not= v_1.$ Since $A \subseteq {\rm cl}\, U,$ this implies $f_w(U)\cap f_v(A)=\emptyset .$ So the boundary set $A\cap f^{-1}_wf_v(A)$ cannot contain a point $x\in U.$ The set $U$ is disjoint to each boundary set, and also to the closure of their union. The last assertion now follows from the fact that $U$ for the OSC can be chosen so that $A\cap U\not=\emptyset$ \cite{Sch}. 

b) The central open set $U$ essentially consists of all points $x$ which have smaller distance from $A$ than from $B(A).$ See \cite{BHR} for details. 

c) Suppose that the topological boundary coincides with the dynamical one.  If the topological  interior of $A$ is empty, then the boundary coincides with $A.$ By (a), this contradicts the OSC. \quad  Conversely, if $A$ has non-empty interior, and the OSC holds, the interior can be taken as open set $U.$ The topological boundary of $A$ then is ${\rm cl}\, U\setminus U.$ Moreover, there is a piece $f_w(A)\subset U.$ Its topological boundary must be covered by sets $f_v(A).$ Because of the OSC, the $f_v(A)$ cannot contain interior points of $f_w(A).$ Thus the boundary sets $A\cap f^{-1}_wf_v(A)$ cover the topological boundary of $A.$ Actually, there is only a finite number of these boundary sets, due to the arguments of P.A.P.~Moran \cite{Mo}, M.~Moran \cite{Mo1}, and Schief \cite{Sch}.
\hfill $\Box$ \vspace{1ex}

Attractors with non-empty interior generate tilings  \cite{Baake2013}, as can be derived from  the above proof by iterating $f_w^{-1}.$
Boundaries of self-similar and self-affine tiles $A$ have been considered in numerous papers. It was discovered that the boundaries have a self-similar structure themselves, represented by a system of equations of type \eqref{hut}, a so-called graph-directed construction \cite{MW}.  This is clear for a $2\times 2$ square, and not difficult for the twindragon \cite{Gi86}.  For the L\'{e}vy dragon with boundary dimension 1.93.., a lot of work was required \cite{DKV,SW}. All these examples were periodic tiles, and there is a unique subdivision of $B(A)$ into boundary sets.\vspace{1ex}

\begin{figure}[h!b]
\begin{center}
\includegraphics[width=0.42\textwidth]{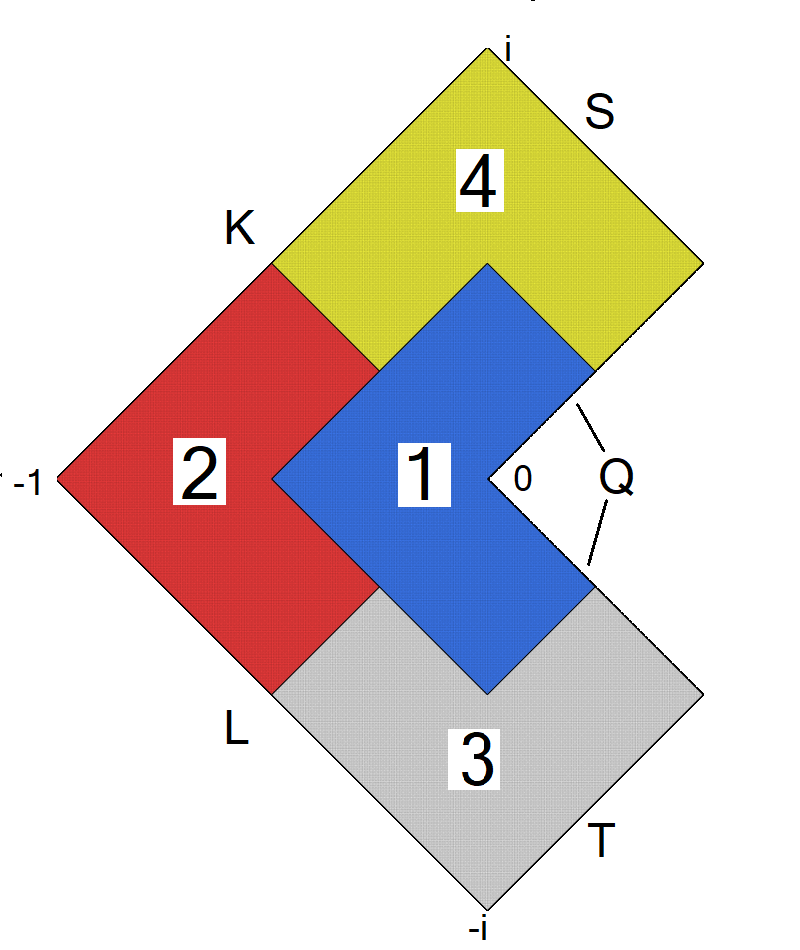}  \quad 
\includegraphics[width=0.42\textwidth]{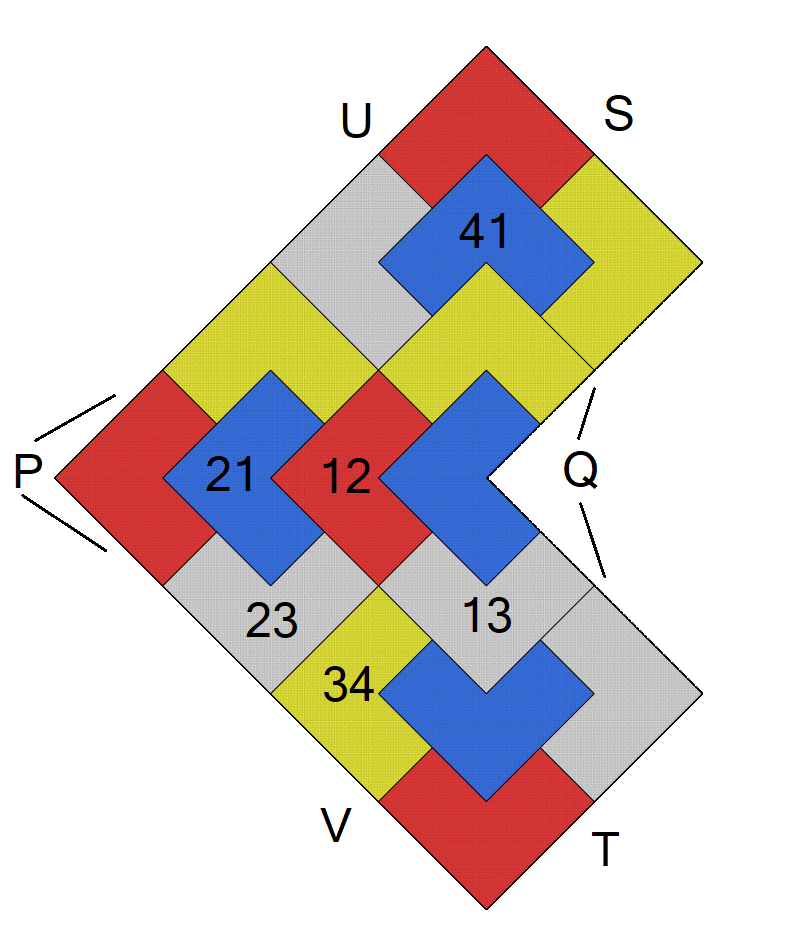}
\end{center}
\caption{The chair, generated by $f_1(z)=\frac{z}{2},\ f_2(z)=\frac12 (z-1),\  f_3(z)=\frac{i}{2} (z-1),$ and $f_4(z)=\frac{-i}{2} (z-1).$  Left: Level 1, with five boundary sets indicated. Right: Level 2, with three other boundary sets.}
\label{chair}      
\end{figure}

{\bf Boundaries of the chair tile.}
For a non-periodic self-similar tile, the situation is more complicated. As a most simple example, we discuss the chair \cite[Chapter 11]{GS} where the boundary consists of line segments.  (For fractal boundaries see \cite[Section 1.7]{EFG3}.)  To find all boundary sets, we need only consider the $A_v\cap A_w$ with $v_1\not= w_1.$ That is, the intersections of first level pieces on the left of Figure \ref{chair}, and their parts on higher levels. 

On the first level, the intersection $D=A_2\cap A_3$  leads to the boundary line segments $T=f_2^{-1}(D)$ and $S=f_3^{-1}(D).$  The boundary sets $Q=f_2^{-1}(A_1\cap A_2)$ and $P=f_1^{-1}(A_1\cap A_2)$ are angles, composed of two segments. So are $M=V\cup T=f_1^{-1}(A_1\cap A_3)$ and  $N=U\cup S=f_1^{-1}(A_1\cap A_4).$

Since $D=A_{23}\cap A_{34},$ we get $L=f_{23}^{-1}(D)$ and  $K=f_{34}^{-1}(D)$ on the second level.  Similarly we have $U=f_{13}^{-1}(A_{13}\cap A_{34})$ and $V=f_{13}^{-1}(A_{13}\cap A_{33}).$ 

Another method goes from one level to the next by subdivision of already collected boundary sets, using the self-similarity of $A.$  Inspection of Figure \ref{chair} gives the set equations
\[ T=f_3(L), \ S=f_4(K), \ Q=f_1(Q)\cup f_3(T)\cup f_4(S)\ . \]
Here $L=f_3^{-1}(T)=f_3^{-1}f_2^{-1}(D)$ and $K=f_4^{-1}(S)$ are new boundary sets on level 2.  The equation for $Q$ involves only the old sets. From $P, M,$ or $N$ we get the new sets $U,V:$
\[ P=f_2(P)\cup f_2(U)\cup f_2(V), \  M=f_3(P\cup  U\cup V), \  N=f_4(P\cup  U\cup V)\, . \]
Finally, the subdivision of $V, U, L,$ and $K$ leads back to $K$ and $L.$
\[ V=f_3(K), \ U=f_4(L), \ L=f_2(L)\cup f_3(K), \  K=f_2(K)\cup f_4(L)\, . \]
Thus on level 3 there are no new sets. This shows that the list is complete. 
The equations can be subsumed by a directed graph, with vertices representing the boundary sets and with labels $i\in \{ 1,...,m\}$ at the edges. The last equation, for instance, is represented by two edges  $K\to K$ and $K\to L$ with labels 2 and 4, respectively.  This is why Mauldin and Williams called this a graph-directed construction \cite{MW}.  Such a labelled graph is also called a finite automaton. For fractal boundaries,  the equations determine the Hausdorff dimension of the boundary.  However, sets like $U$ and $V$ are not explicitly defined by the equations. It is only known that the equations have a unique solution within compact sets. \vspace{1ex}

{\bf The neighbor graph concept.}
We shall consider another automaton which is called the neighbor graph of the IFS. The decisive step is to replace a boundary  $H=A\cap f_w^{-1}f_v(A)$ by the mapping $h=f_w^{-1}f_v$ which maps $A$ to the neighbor set $h(A)$ which shares the boundary $H$ with $A.$ The map $h$ is called neighbor map. The intuition is that we want to extend the fractal construction outwards with a set of puzzle pieces, all congruent to $A.$ The neighbor maps prescribe the position of the next piece in such a way that boundary sets meet each other.

The replacement of boundary sets by neighbor maps is a tremendous simplification. We can calculate with maps, not with boundary sets. The determination of the neighbor graph can be done by computer, without any inspection of figures. Only in this way we can begin to study more complicated examples \cite{EFG3,EFG5,TZ20}. Moreover, the neighbor graph completely describes the topology of the attractor $A$ \cite{EFG4}. \vspace{1ex}

\begin{figure}[h!b]
\begin{center}
$N^1$\!\includegraphics[width=0.216\textwidth]{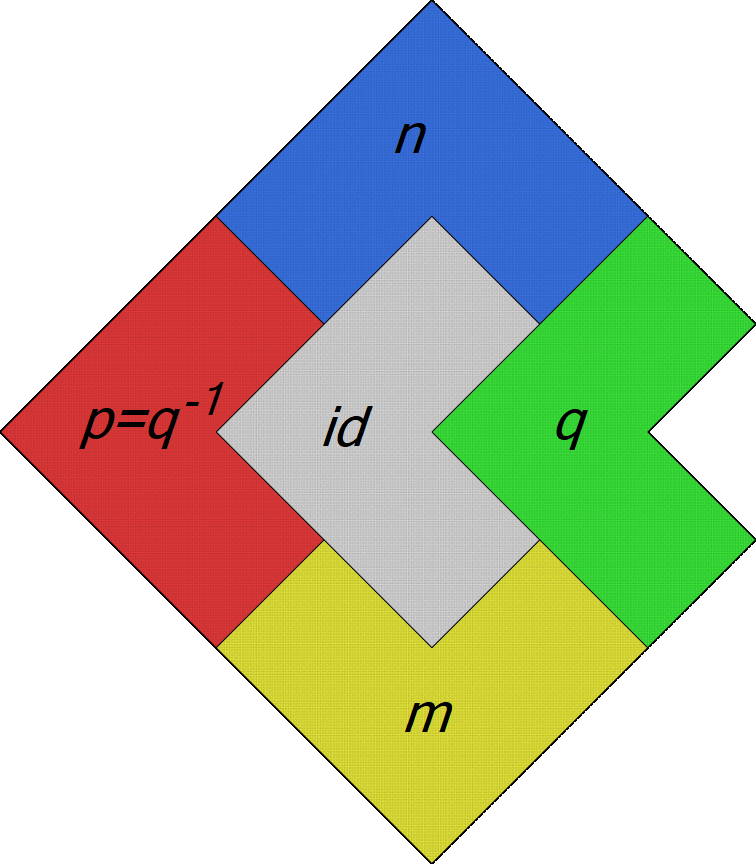}  \ 
$N^{4a}$\includegraphics[width=0.29\textwidth]{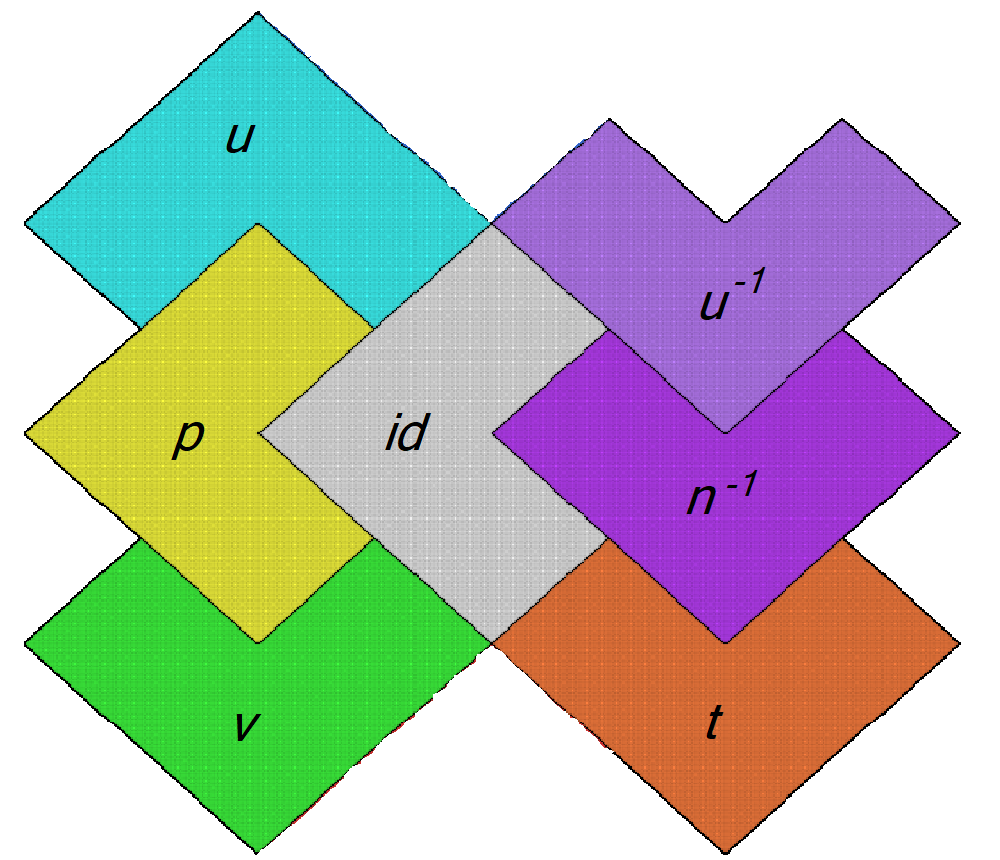} \quad 
$N^{3b}$\includegraphics[width=0.28\textwidth]{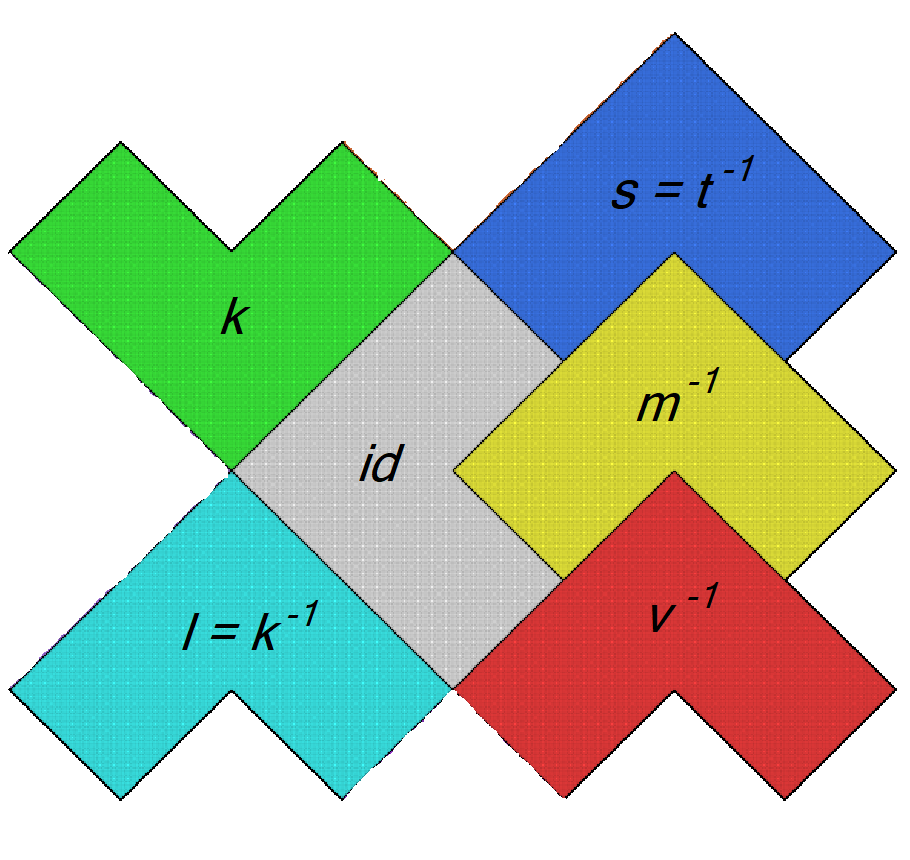} 
\end{center}
\caption{The positions of neighboring tiles of the chair, denoted by the neighbor maps. The boundary set $Q$ corresponds to three neighbor maps $q=p^{-1}, m^{-1},$ and $n^{-1}.$ These neighborhoods appear in every tiling and will be called $N^1,N^{4a}$ and $N^{3b}$ in Section \ref{inte}.}
\label{chairnb}      
\end{figure}

{\bf The neighbor graph of the chair tile.}
For the chair example, we chose the coordinate system in the complex plane so that the IFS is $f_1(z)=\frac{z}{2},\ f_2(z)=\frac12 (z-1),\  f_3(z)=\frac{i}{2} (z-1),$ and $f_4(z)=\frac{-i}{2} (z-1).$ Neighbor maps are denoted by small letters $p,q,s,...$ corresponding to the boundary sets $P,Q,S,... .$ They all have the form $h(z)=az+b$ with $a\in\{ 1,i,-i\}$ and $b\in \{\pm 1,\pm i, \pm 1\pm i \}.$ In Figure \ref{chairnb} the neighbor set $h(A)$ is denoted by $h,$ and the reference tile $A$ by $id,$ the identity map. Note that the inverse of a neighbor map $h$ must be a neighbor map, with $h(A)$ and $A$ interchanged. 

The calculation of neighbor maps, described in detail in \cite{EFG3,EFG5,EFG1}, is similar to the calculation of boundaries above.  On the first level we determine the $f_i^{-1}f_j.$ 
\[ p=f_1^{-1}f_2=z-1, \ q=f_2^{-1}f_1=z+1, \ t=f_2^{-1}f_3=iz+1-i, \ s=f_3^{-1}f_2=-iz+1+i ,  \]
\[ m=f_1^{-1}f_3=iz-i, \, n=f_1^{-1}f_4=\! -iz+i, \, m^{-1}\! =f_3^{-1}f_1=\! -iz+1,
\, n^{-1}\! =f_4^{-1}f_1=\! iz+1.\]
The neighbor maps $q=p^{-1}, m^{-1}$ and $n^{-1}$ all describe the boundary set $Q,$ but with different neighbor positions. See Figure \ref{chairnb}.  The recursive calculation of further maps $h'$ from the $h$ in our list proceeds with the formula 
\begin{equation}h'=f_i^{-1}hf_j \ .
\tag{4}\label{h'h}\end{equation}
Similarly to $L=f_3^{-1}(T)$ and $K=f_4^{-1}(S)$ we get
\[ l =f_3^{-1}tf_4= f_{23}^{-1}f_{34}=-iz-1-i \quad \mbox{ and } \quad  k =f_4^{-1}sf_3= f_{34}^{-1}f_{23}=iz-1+i \, .\]
Finally, $u=f_3^{-1}mf_4$ and $v=f_3^{-1}mf_3.$ The maps $u^{-1}$ and $v^{-1}$ belong to the boundary sets $S$ and $T,$ respectively.  Altogether we found 14 neighbor maps for 10 boundary sets.  For simplicity, we have restricted ourselves to maps $h=f_w^{-1}f_v$ for which $A\cap h(A)$ consists of more than one point. There are five other maps where the intersection is a singleton.

\begin{figure}[t]
\begin{center}
\includegraphics[width=0.748\textwidth]{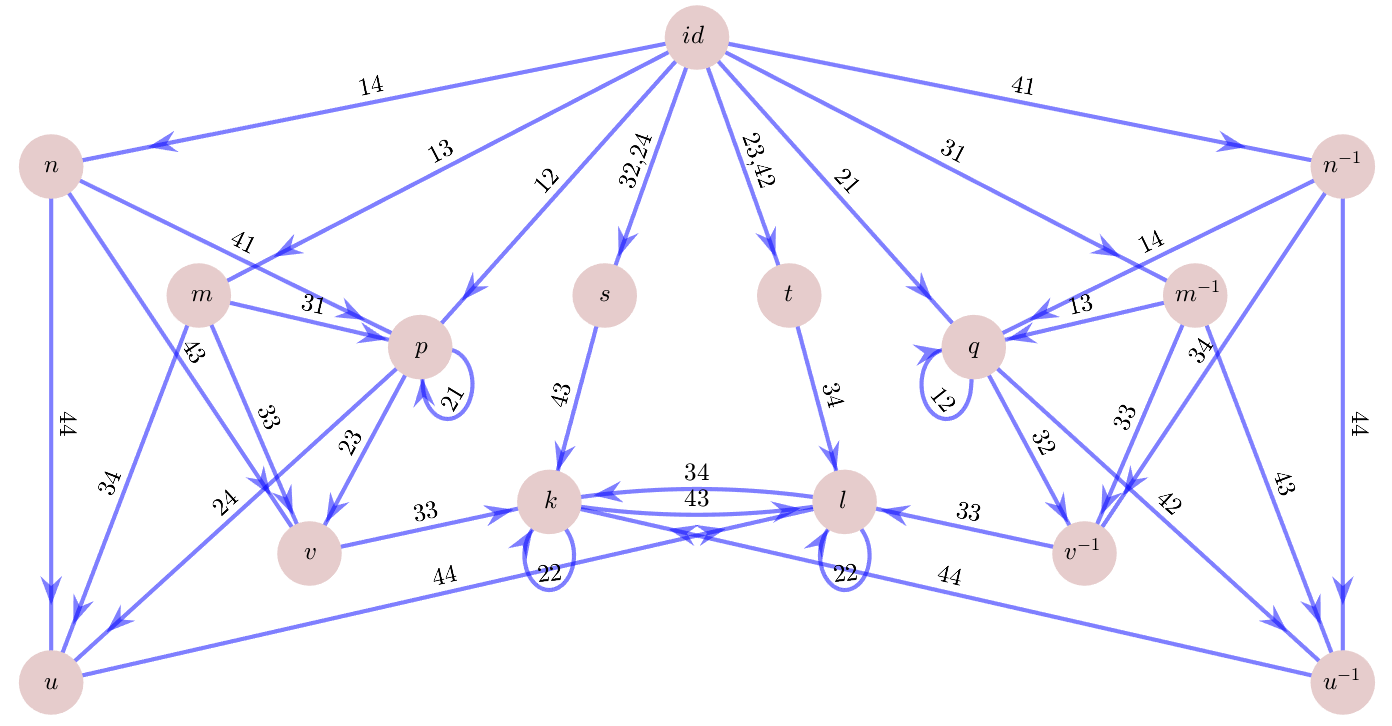}  
\end{center}
\caption{The neighbor graph for the chair tile, with 15 vertices and 38 edges. }
\label{chairng}      
\end{figure}

The neighbor graph has as vertices $id$ and all neighbor maps $h.$ From $id$ there are edges to the first level neighbor maps $h=f_i^{-1}f_j,$ with corresponding label $i\, j.$ When the relation \eqref{h'h} holds for two neighbor maps $h,h',$ there will be an edge from $h$ to $h'$ with label $i\, j.$ For the chair, only some of the relations have been mentioned so far. Figure \ref{chairng} shows the complete neighbor graph. For computer work, a list of edges $(h,h',i,j)$ is more appropriate -- a $4\times n$ matrix, where $n$ is the number of edges. By definition, there are loops at $id$ with label $i\, i$ for each $i.$ They are not  shown here. Instead of drawing multiple edges, we write multiple labels at an edge.

There is an obvious correspondence between neighbor maps and their inverses \cite{EFG4}, given in Figure \ref{chairng} by the reflection at the vertical line through $id.$ We shall use the graph in the next section.

\begin{figure}[t]
\begin{center}
\includegraphics[width=0.3\textwidth]{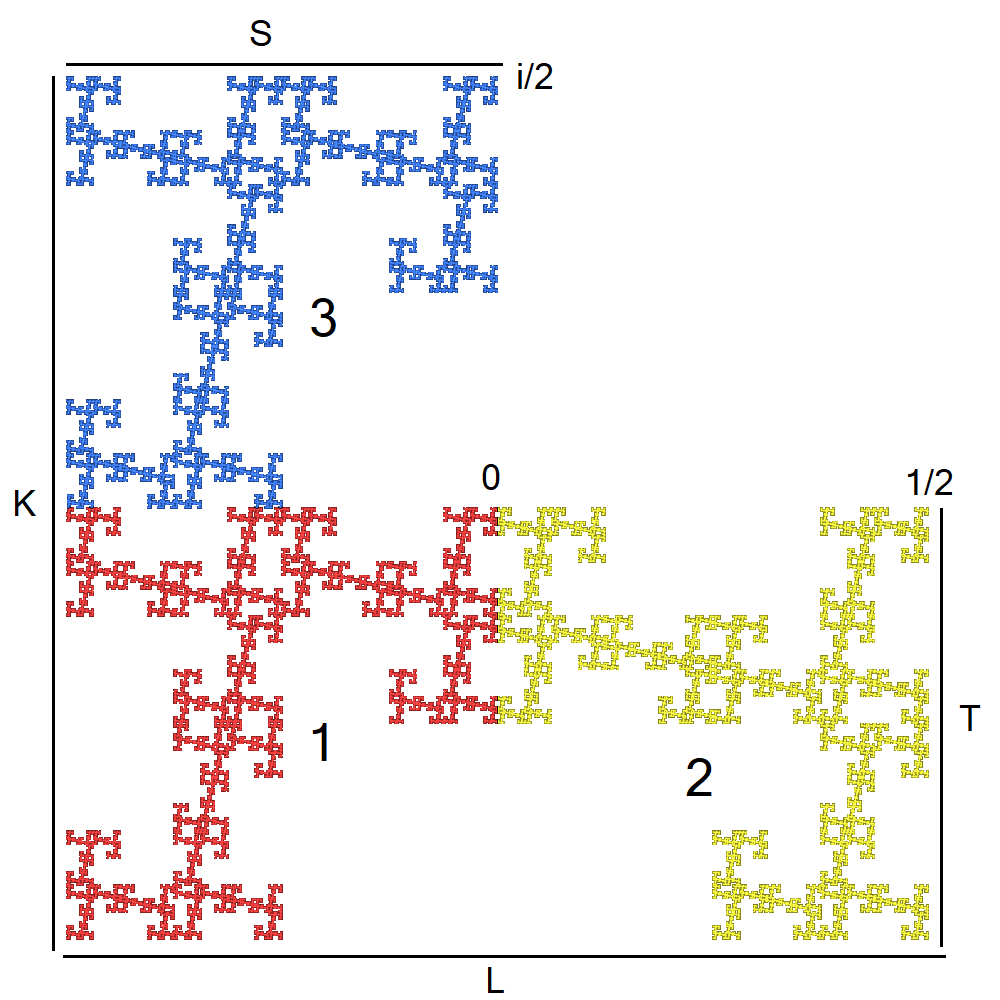}  \quad
\includegraphics[width=0.3\textwidth]{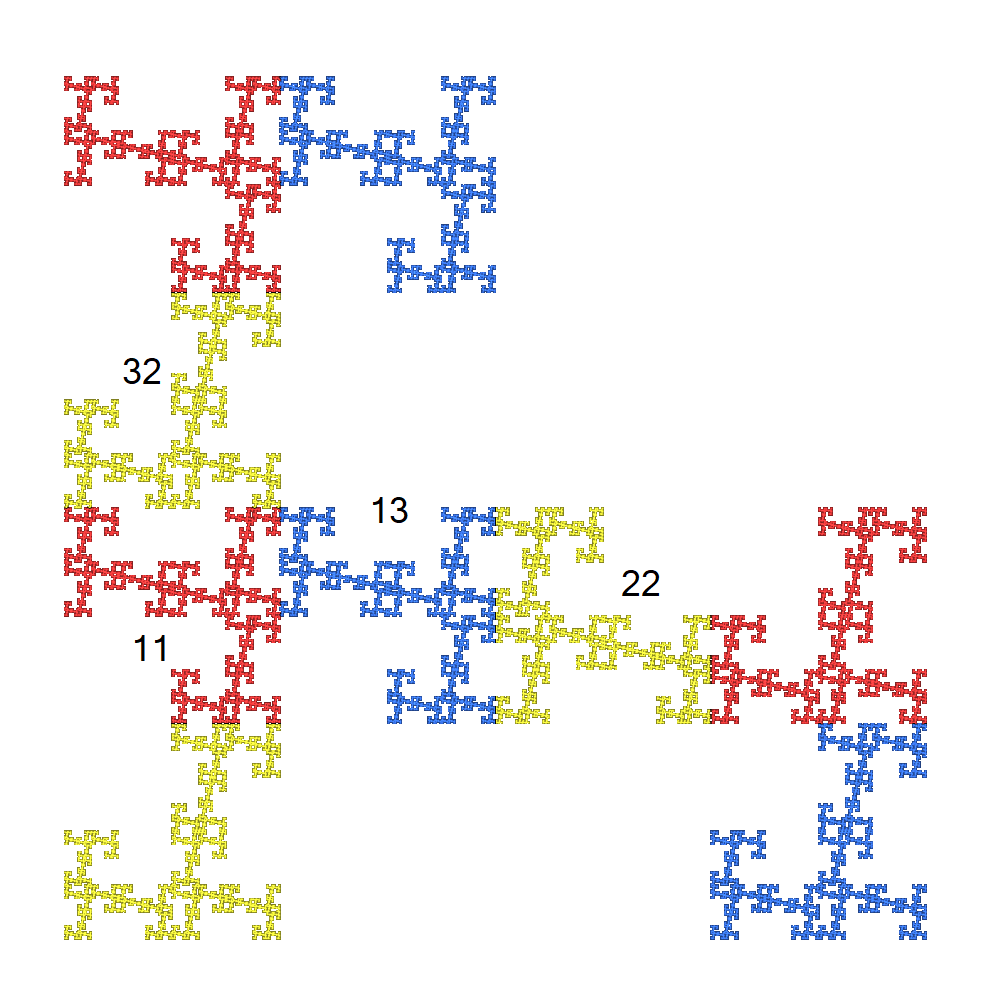} \quad 
\includegraphics[width=0.3\textwidth]{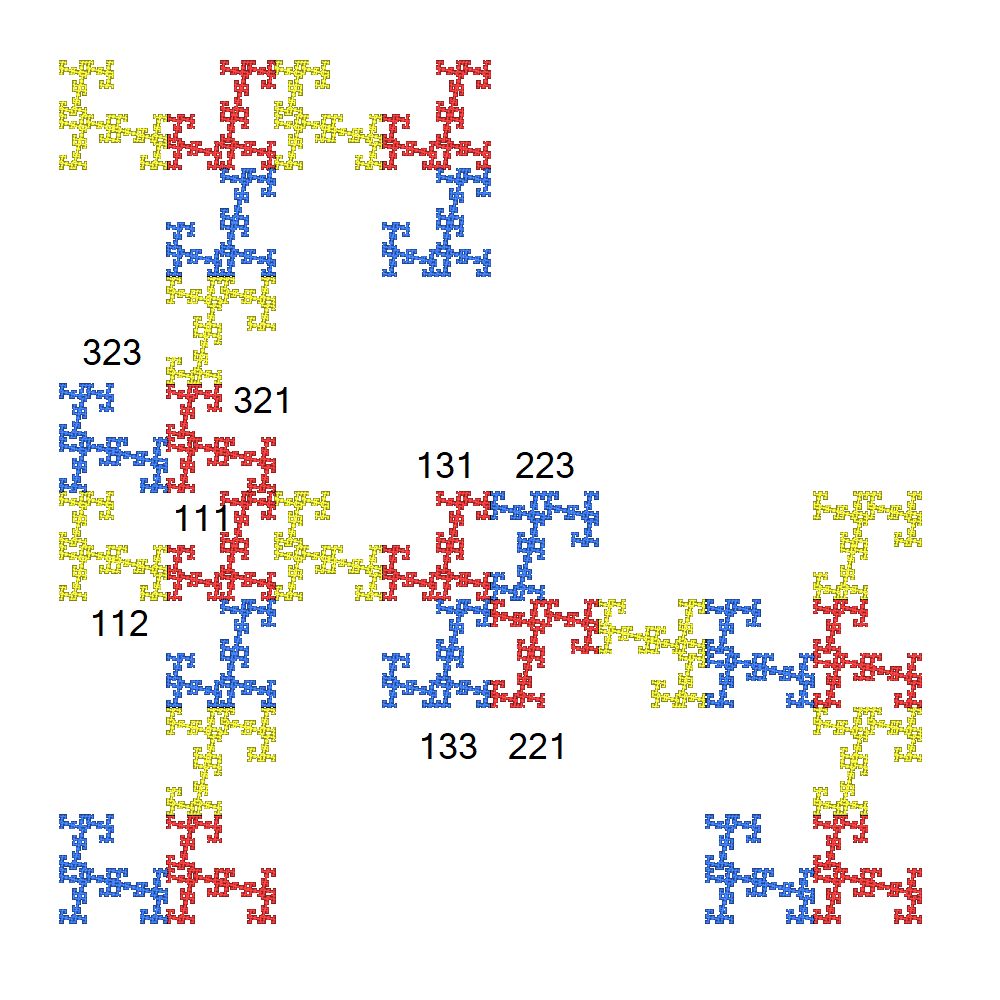} 
\end{center}
\caption{A self-similar set with Cantor boundary sets of two different dimensions. }
\label{figH}      
\end{figure}

{\bf A fractal square.}
In contrast to tiles where topological and dynamical boundaries coincide, we now consider a proper fractal in the plane.  For simplicity, we take Figure \ref{figH} where a square can be taken as open set for the OSC. The boundary sets are then contained in the sides of the square. They are denoted $S, T, K, L$ as for the chair. Now they are Cantor sets. Nevertheless, boundary equations can be established as for the chair:
\[ S=f_3(K), \ T=f_2(K), \ L=f_1(T)\cup f_2(S), \ K= f_1(L)\cup f_3(L)\, .\]
Substituting for $L,T,S$ we obtain $K=f_{112}(K)\cup f_{123}(K)\cup f_{312}(K)\cup f_{323}(K)$ which can be seen in the figure. So $K$ is a self-similar set with four mappings of contraction factor $1/8.$ Its dimension is $\log 4/\log 8 =2/3.$ 

However, these were only the boundary sets associated with $f_1(A)\cap f_2(A).$ The two boundary sets from $f_1(A)\cap f_3(A)=f_{32}(K)\cap f_{11}(L)$ have smaller dimension.  They will be determined from the neighbor graph.

\begin{figure}[b]
\begin{center}
\includegraphics[width=0.7\textwidth]{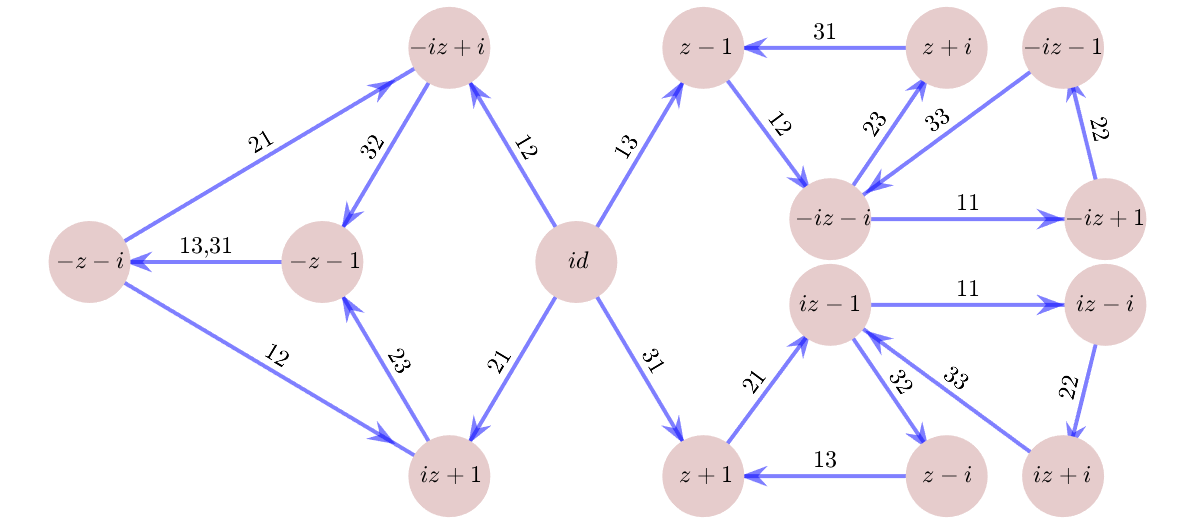}
\end{center}  
\caption{The neighbor graph for the fractal square in Figure \ref{figH}.  }
\label{figHng}     
\end{figure}

The coordinate system, indicated by the unit points $0, 1/2, i/2,$ was chosen so that the IFS 
is $f_1(z)=\frac{-iz}{2}-\frac{1+i}{4}$ (red), $f_2(z)=\frac{-z}{2}+\frac{1-i}{4}$  (yellow), and $f_3(z)=\frac{-iz}{2}-\frac{1-i}{4}$  (blue).  With this choice, all neighbor maps have the form  $h(z)=az+b$ with $a,b\in\{ 1,i,-1,-i\}.$ First the square is rotated around zero and then translated to a neighboring square. 
Two of these 16 maps, namely $-z+1$ and $-z+i,$ do not appear in the construction.
The neighbor graph in Figure \ref{figHng} can be established by inspection of Figure \ref{figH} and/or by calculation. We only consider neighbors which intersect in more than one point.

The correspondence between neighbor maps and their inverses is given here by reflection at the horizontal symmetry axis. The component on the left corresponds to the boundary sets calculated above, for instance $k(z)=f_{13}^{-1}f_{22}=-z-1$ belongs to $K.$ The four cycles from $k$ to $k$ have first label words 123, 323, 112, 312 as in the equation above. On the right, the mapping  $f_{32}^{-1} f_{11}= iz-1$ is the starting point for two cycles to itself, with first label words 312 and 123. This means that 
$B=A\cap f_{32}^{-1} f_{11}(A)$ is a self-similar set $B=f_{312}(B)\cup f_{123}(B)$ with two maps with contraction factor $1/8.$  The Hausdorff dimension is $\log 2/\log 8 =1/3.$ Note that $B\subset K$ while the boundary set of the inverse map $f_{11}^{-1} f_{32}= -iz-i$ is a subset of $L.$ The system of boundary sets in Figures \ref{figA}, \ref{figB}, \ref{figC} is much more complicated.

\section{Dynamical interior and neighborhood graph}  \label{inte}
{\bf The dynamical interior as a union of pieces.}
For a self-similar set $A$ with open set condition, the dynamical interior $I(A)$ is defined as the complement of the dynamical boundary in $A.$  We give a description of $I(A)$ in terms of the pieces of $A.$ 

For each piece $A_w,\ w=w_1...w_k,$ the \emph{neighborhood} consists of all pieces $A_v$ of the same size which intersect $A_w.$ If contraction factors $r_1,...,r_m$ are equal, equal size means that $v=v_1...v_{k'}$ has length $k'=k.$ If the $r_i$ are different, $A_v$ will have ``the same size'' as $A_w$ if $r_v=r_{v_1}\cdots r_{v_{k'}}\ge r_w$ and $r_v<r_{w_1}\cdots r_{w_{k-1}}.$ That is, the neighbor piece is allowed to be a little larger than the reference piece.  
With this convention, the neighborhood of $A_w$ is unique.
(It would have been possible to define neighbors and the neighbor graph for different contraction factors, too.  However, it would be counterintuitive to work with a non-symmetric neighbor relation. For the next theorem, there are no such problems.)

We say that the neighborhood $N=\{ A_{v^1},..., A_{v^n}\}$ of $A_w$ is maximal if there is no piece $A_u$ which has the same neighbors, up to similarity (that is, the neighbors
$f_uf_w^{-1}(A')$ with $A'\in N$), and at least one more neighbor. 
 
\begin{Theorem}\label{intepie}
For a self-similar set in $\RR^d$ with OSC,  the following is true.
\begin{enumerate}
\item[a)]  The dynamical interior $I(A)$ is a countable union of pieces.
\item[b)]  All pieces with maximal neighborhood are contained in $I(A).$
\item[c)]  For any two pieces $A_v, A_w$ inside $I(A),$ there is a subpiece $A_u\subset A_v$  such that the neighborhoods of $A_u$ and $A_w$ are geometrically similar.
\end{enumerate}
\end{Theorem}

{\it Proof. } 
a) For a point $x\in I(A),$ let $d$ be the distance from $x$ to the dynamical boundary.  There is an integer $k$ such that all pieces of level $k$ have diameter smaller than $d.$ The union of pieces of level $k$ containing $x$ is a subset of $I(A).$ Thus $I(A)$ is a union of pieces, and there are only countably many pieces.

b) Let $y$ belong to a piece $A_w$ with maximal neighborhood $N.$ We show that $y$ cannot belong to the dynamical boundary of $A.$ If $y$ would be in a boundary set $f_u^{-1}f_v(A),$ then $A_u$ would intersect the pieces $f_v(A)$ and $f_uf_w^{-1}(A')$ with $A'\in N.$ This would contradict the assumption that $N$ is maximal.  Because of the OSC, there can be only finitely many boundary sets which intersect $N$ \cite{Sch,BSS}. Thus $y$ cannot be an accumulation point of a sequence of boundary sets, and thus belongs to $I(A).$

c) Take $u=vw.$ Then $A_u\subseteq A_v\subseteq I(A),$ and $A_u$ has all the neighbors of $A_w$ since $A_{vw}\cap A_{vw'}=f_v(A_w\cap A_{w'}).$  There are no more neighbors since $f_v^{-1}(A_u)=A_w$ does not intersect the boundary of $A.$
\hfill $\Box$ \vspace{1ex}

Note that there can be neighborhoods which occur for pieces inside $I(A)$ as well as for pieces which intersect $B(A).$  As an example, take the well known Sierpi\'nski carpet, with the central hole in a $3\times 3$ square.  The boundary consists of the four sides of the big square. Pieces inside $I(A)$ appear only at the second level. At first level, the $A_i$ at the four corners have two neighbors. Such corner neighborhoods can never occur inside $I(A).$ The four other pieces $A_i$ also have two neighbors, either in a row or in a column. (We neglect neighbors with one common point.)  Lots of smaller copies of these neighborhoods appear inside $I(A),$ between two holes of different size.\vspace{1ex}

{\bf Finding an interior piece.} Now we assume that the IFS has equal contraction factors, fulfils the OSC and is of finite type. Thus there are $n$ neighbor maps for some integer $n.$ Each neighborhood is a collection of neighbors, or neighbor maps. So there can be at most $2^n$ different neighborhoods.  We shall present a simple algorithm which determines all neighborhoods of interior pieces directly from the neighbor graph.

In a first step, we have to find one piece inside $I(A).$ For this purpose, we have to study the directed edge paths starting in any vertex of the neighbor graph except $id.$ Each  vertex corresponds to a boundary set, and the first labels on the edge paths describe the words $w$ of pieces $A_w$ which intersect this boundary set. (The second labels describe the boundary from the other side. They coincide with the first labels of the `inverse' path.) See for example \cite{BM09,EFG3,EFG4}.  We have to find a word which does \emph{not} occur as labelling of an edge path starting in any vertex except $id.$
\vspace{1ex}

For the chair example, consider Figure \ref{chairng} and remove all the edges starting in $id.$ In vertex $s$ starts only one edge to $k,$ with first label 4. That means the boundary set $S$ is contained in $A_4,$ as seen in Figure \ref{chair}. From $k$ there are two edges with first label 2 and 3. That is, $K$ intersects $A_2$ and $A_3.$ If we combine these edges with the edge from $s$ to $k,$ we see that $S$ intersects $A_{42}$ and $A_{43}.$ Thus $A_4, A_2, A_3, A_{42},$ and $A_{43}$ are not interior pieces.  The label 1 occurs only at edges from $q, m^{-1},$ and $n^{-1}.$ There are no edges which end at $m^{-1}$ or $n^{-1},$ and all edges which end at $q$ have first label 1.  Thus the path $w1$ occurs only as an edge path if $w=11...1.$  So the pieces $A_{21}, A_{31}, A_{41}, A_{121},...$ and of course their subpieces all belong to the interior $I(A).$ We could have seen this in Figure \ref{chair}, but our algorithm must not rely on geometric arguments.

For the fractal square, look at Figure \ref{figHng}. There are edges with label 2 from $iz+1$ and five other vertices. But these edges do not connect any two of these vertices. So the label 22 does not occur on a path of two edges. Figure \ref{figH} shows that indeed $A_{22}$ and $A_{13}$ are the two interior pieces on second level.\vspace{2ex}

{\bf Construction of the neighborhood graph.}
Once we have an initial interior piece $A_v$ with neighborhood $N_v,$ we can apply Theorem \ref{intepie}c. We need only study the subpieces of $A_v$ in different levels, in order to get all neighborhoods of interior pieces, up to similarity. Let $N_{vi}$ denote the neighborhood of subpiece $A_{vi}$  for $i=1,...,m.$ Both $N_v$ and $N_{vi}$ are considered as sets of neighbor maps, that is, as sets of vertices of the neighbor graph.

The question now is how to obtain $N_{vi}$ from $N_v.$ First, $A_{vi}$ gets all the neighbors of the first level piece $A_i$ from its sister pieces inside $A_v.$ That is, $N_{vi}\supseteq N_i.$ Next, $A_{vi}$ inherits some neighbor subpieces from $A_v.$ Let $h=f_v^{-1}f_w\in N_v$ denote a neighbor $A_w$ of $A_v,$ identified with a vertex $s$ of the neighbor graph. For each edge with double label $i\, j$ starting in $s$, the terminal vertex $s'$ represents the neighbor map $h'=f_{vi}^{-1}f_{wj}$ which corresponds to the neighbor $A_{wj}$ of $A_{vi}.$ \vspace{1ex} 

The formal description shows that we do not calculate with maps. Let $G=(V,E)$
be the neighbor graph.  For vertex sets $M\subset V,$ we define the $i$-successor set as
\[  S(M,i) =\{ v'\in V\,|\mbox{ there is a }v\in M\mbox{ and an edge with first label }
i\mbox{ from }v\mbox{ to }v'\} \, .\] 
On the first level, we determine $N_i=S(\{ id\}, i)$ for $i=1,...,m.$  Then we start with an interior piece $A_w$ and its neighborhood $N_w,$ as determined above. We calculate
\begin{equation}\label{nbh}
N_{wi}= N_i\cup S(N_w,i) \ \mbox{ for } i=1,...,m\, . \tag{5}
\end{equation}
Next, the words $wi$ are considered as new $w$ and the formula \eqref{nbh} is applied again.

In this way, a neighborhood graph $G_N=(V_N,E_N)$ is recursively constructed. The vertex set $V_N$ consists of all $N\subset V$ which represent neighborhoods of interior pieces. We start with $V_0=\{ N_w\}$ where $N_w$ is the neighborhood of a first interior piece. Then we check whether the $N_{wi}$ are already in $V_0.$ This would be the case for a periodic self-similar tile. Those $N_{wi}$ which are not in $V_0$ are added to $V_0,$ and the result is $V_1.$ Moreover, an edge with label $i$ is drawn from $N_w$ to $N_{wi},$ for $i=1,...,m.$ Afterwards, the vertex $N_w$ is marked as processed and completed.

For a periodic self-similar tile, where there is only one interior neighborhood, this would finish the construction with $V_1=V_0$ and with a neighborhood graph with one vertex and $m$ loops. If $V_1\not= V_0,$ we take the first unmarked vertex in $V_1$ as new $w$ and apply the recursion \eqref{nbh} again, obtaining $V_2.$ And so on. 

In each recursion step, we get at most $m$ new neighborhoods $N_w,$ and we process and mark one neighborhood $N_v.$ Since there are only finitely many possible neighborhoods, the algorithm finishes at some step $k$ with $V_k=V_{k-1}.$\vspace{1ex}

This procedure is easy to program. It runs fast, even in the case of large numbers of interior neighborhoods.  It is convenient to assign to each $N_v$ an integer. Then
the neighborhood graph $G_N$ with $K$ vertices $N_v$ is represented by a $K\times m$ matrix $M=(m_{ki})$ where $k=1,...,K$ denotes the number of a neighborhood $N_v,$ and $m_{ki}$ is the integer assigned to the successor neighborhood $N_{vi},$ for $i=1,...,m.$ The definition of the vertices $k$ as subsets $N_v\subset V$ is given as a separate matrix.
\vspace{1ex}

{\bf Frequency of neighborhoods.} Pieces with different neighborhoods have different geometry, as noted by Ngai and Wang \cite{NW} for overlapping constructions. The neighborhood graph describes a kind of branching process or substitution. The substitution matrix $S$ of the neighborhoods $N^1,...,N^K$ is the following $K\times K$ matrix.
\begin{equation}
S=(s_{kp}) \mbox{ where } s_{kp} \mbox{ counts the pieces } i \mbox{ with }
N_v=N^k \mbox{ and } N_{vi}=N^p\, .  \tag{6} \label{subst}
\end{equation}
The definition does not depend on $v,$ but requires that $A_v$ is in $I(A).$ In our two examples $S$ contains only zeros and ones. There is a theory for substitution matrices \cite{Q}. The largest eigenvalue $\sigma$ is a positive real and determines the growth factor of the branching process. Ngai and Wang showed that their overlapping attractor has Hausdorff dimension $\log\sigma /-\log r.$ This also holds for our non-overlapping attractors where $\sigma =m$ and the dimension formula is well-known.  The corresponding eigenvector $p=(p_1,...,p_K)$  with $\sum p_k=1$ gives the relative frequencies of neighborhood types in a fractal tiling \cite{GS,NW,Q}.\vspace{1ex}

{\bf Neighborhoods of the chair example.}
For the chair tile, three neighborhoods are shown in Figure \ref{chairnb}.  As initial interior pieces we took $A_{i1}$ for $i=2,3,$ or $4,$ with neighborhood $N^1=\{ m,n,p,q\},$ shown on the left of Figure \ref{chairnb}. We use a superscript to distinguish these special patterns from $N_w$ introduced above. Thus $N^1=N_w$ for $w\in\{ 21,31,41\}$ as shown above and in Figure \ref{chair}.  We determine $N_{w1}$ by \eqref{nbh}. The set $N_1$ of 1-successors of $id$ in Figure \ref{chairng} is $\{ m,n,p\} ,$ and these are just the neighboring pieces of $A_1$ in the first level.  We have to add the 1-successors of $N^1.$ This is only the vertex $q$ because of its loop. Thus $N_{w1}=N^1.$ The neighborhood graph starts with a loop with label 1 from $N^1$ to itself. Similarly,  $N_{w4}$ is the union of $N_4=\{ t, n^{-1}\}$ and of the 4-successors of $N^1,$ that is $\{ p,u,v,u^{-1} \}.$ The result is the neighborhood $N^{4a}$ shown in the middle of Figure \ref{chairnb}. Thus in $G_N$ there is an edge with label 4 from $N^1$ to $N^{4a}.$ Next, an edge with label 3 will lead from $N^{4a}$ to $N^{3b}$ drawn on the right of Figure \ref{chairnb}. The words $w$ are not needed for the calculation, only the sets like $N^{4a},$ and the neighbor graph. The complete neighborhood graph is given in Table \ref{chairnbh}. Besides $N^1$ there are six neighborhoods. They differ in their left-hand parts $a=\{ p,u,v\} ,\ b=\{ k,l\},$ or in the right-hand parts $\{ q,s,t\} ,\ \{ m^{-1},s,v^{-1}\} , $ and $\{ n^{-1},t,u^{-1}\} ,$ which correspond to 2,3,4, respectively, as last letter of $v.$\vspace{1ex}

\begin{table}[h] \label{chairnbh}       
\begin{center}
\begin{tabular}{|l|l|llll|l|}
\hline
&Neighborhood & Successors & & & &Frequency\\
Name&Neighbors&1&2&3&4& \\ \hline
$N^1$ & $m,n,p,q$ & $N^1$ & $N^{2a}$ & $N^{3a}$&$N^{4a}$&$1/4$\\
$N^{2a}$ & $p,q,s,t,u,v$ & $N^1$ & $N^{2a}$ & $N^{3b}$&$N^{4b}$&$1/8$\\
$N^{3a}$ & $m^{-1},p,s,u,v,v^{-1}$ & $N^1$ & $N^{2a}$ & $N^{3b}$&$N^{4b}$& $1/16$\\
$N^{4a}$ & $n^{-1},p,t,u,u^{-1},v$ & $N^1$ & $N^{2a}$ & $N^{3b}$&$N^{4b}$& $1/16$\\
$N^{2b}$ & $k,l,q,s,t$ & $N^1$ & $N^{2b}$ & $N^{3b}$&$N^{4b}$&$1/8$\\
$N^{3b}$ & $k,l,m^{-1},s,v^{-1}$ & $N^1$ & $N^{2b}$ & $N^{3b}$&$N^{4b}$&  $3/16$\\
$N^{4b}$ & $k,l,n^{-1},t,u^{-1}$ & $N^1$ & $N^{2b}$ & $N^{3b}$&$N^{4b}$&  $3/16$\\
\hline
\end{tabular} \end{center}
\caption{The graph of the seven interior neighborhoods of the chair tile.}
\end{table}

We calculate the relative frequencies $p(N)$ of neighborhoods in a large chair tiling, as drawn in \cite[Chapter 11]{GS}.  Since there are equal numbers of tiles which occupy position 1,2,3, or 4 in their supertile, we have $p(N^1)=1/4= p(N^{ia})+ p(N^{ib})$ for $i=2,3,4.$  For the $7\times 7$ substitution matrix corresponding to Table \ref{chairnbh} equation \eqref{subst} gives $p(N^{2a})=p(N^{2b})=1/8$ and $p(N^{3a})=p(N^{4a})=1/16, \ p(N^{3b})=p(N^{4b})=3/16.$ The neighbors $k,l$ appear three times more often than $p,u,v$ in pieces of type 3 and 4 but equally often in pieces of type 2.\vspace{1ex}

\begin{figure}[h]
\begin{center}
\includegraphics[width=0.25\textwidth]{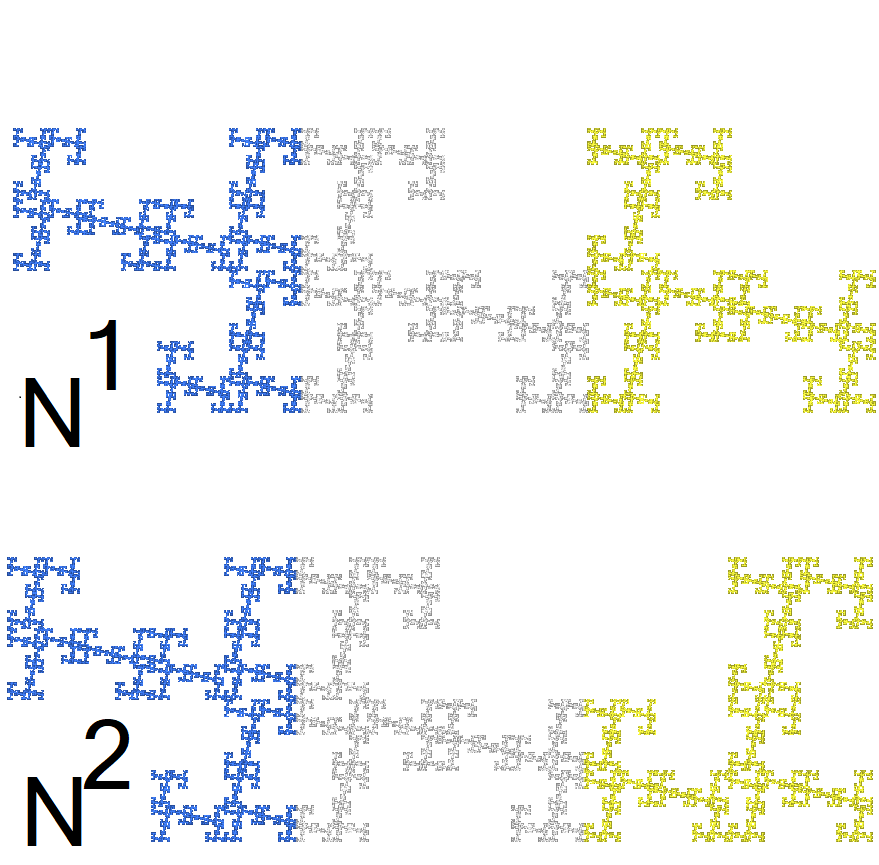}  \qquad
\includegraphics[width=0.25\textwidth]{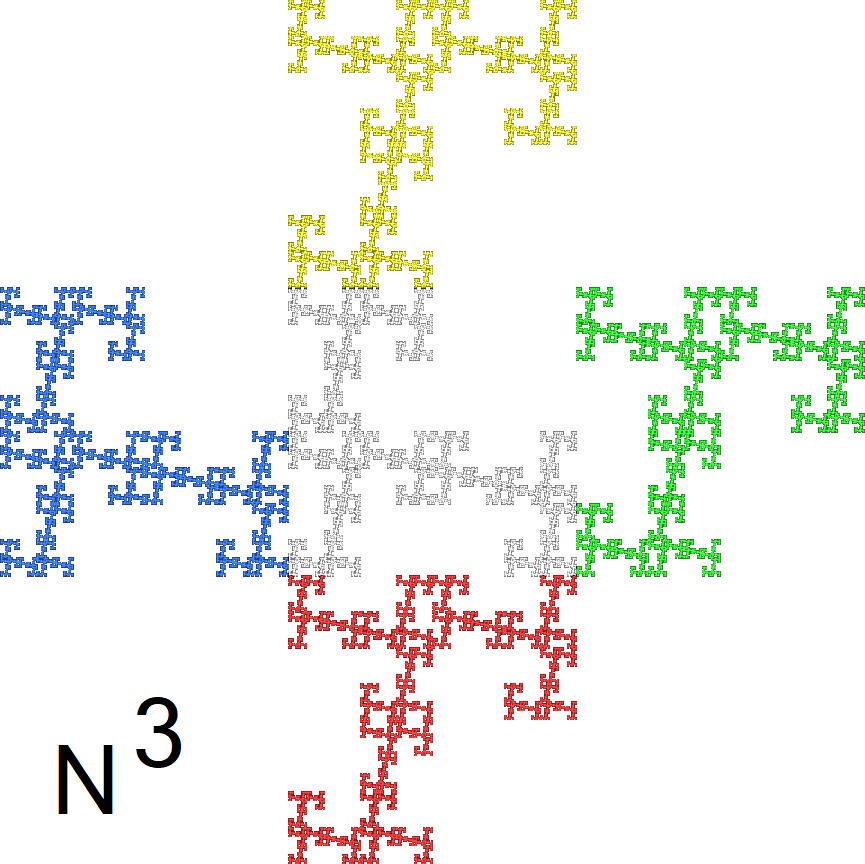} \qquad 
\includegraphics[width=0.25\textwidth]{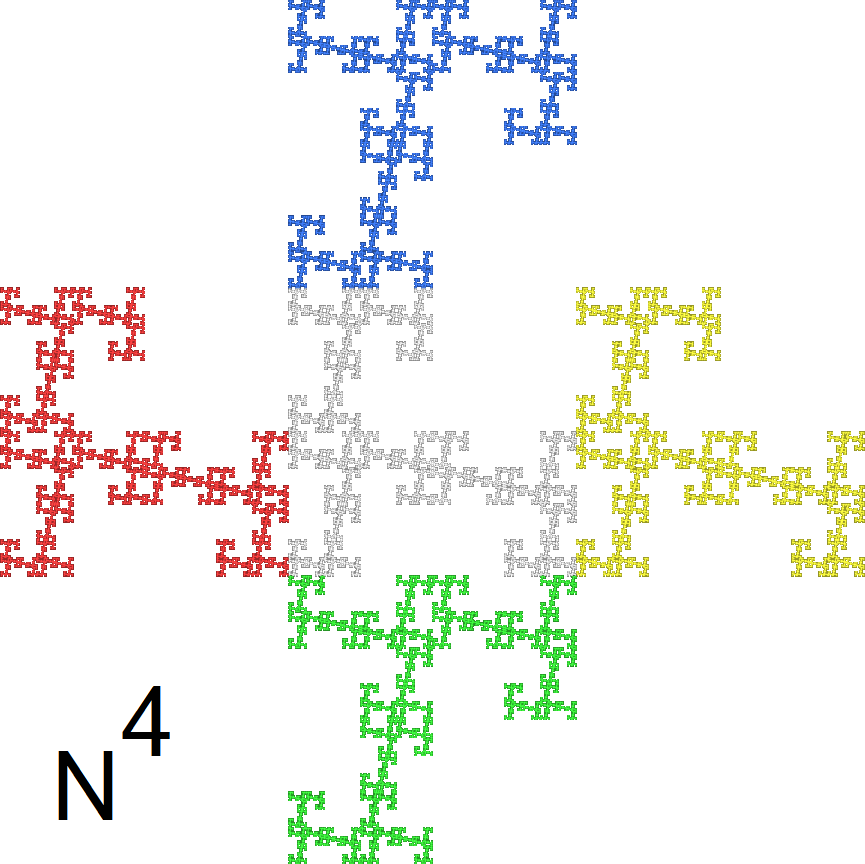} 
\end{center}
\caption{Four of the 30 neighborhoods of the fractal square, with the reference piece in grey. }
\label{figHnbh}      
\end{figure}

{\bf Neighborhoods for the fractal square.}
For self-affine or crystallographic tilings, we have a single neighborhood because of the transitivity of the symmetry group. In case of an aperiodic tile like the chair, there is a small number of different neighborhoods. For fractal examples, the number of neighborhoods is much larger, except for toy models like the Sierpi\'nski gasket. If there are $n$ neighbors, there could be $2^n$ neighborhoods.  For our fractal square, we found 14 neighbor maps $az+b$ with $a,b\in\{ \pm 1,\pm i\}$ in Figure \ref{figHng}.  Since each neighbor fits either on the left, right, top or bottom of the square, the bound $2^{14}$ reduces to $4\cdot 4\cdot 3\cdot 3=144.$ The calculation of the neighborhood graph gave 30 neighborhoods.  Instead of a table, we show four neighborhoods in Figure \ref{figHnbh}.  As for the chair, we determined the frequencies of the neighborhoods with the substitution matrix \eqref{subst}. This allows a brief statistical analysis of the dynamical interior. \vspace{1ex} 

Figure \ref{figHnbh} shows on the left the two most frequent neighborhoods $N^1$ and $N^2,$ each with a frequency of 12 percent. Altogether there are eleven neighborhoods with two neighbors which account for 56 percent of all neighbor patterns. Thirteen patterns with an overall frequency of 28 percent have three neighbors. Four neighborhoods with total frequency of 7 percent have the maximum number of four neighbors. The neighborhood $N^3$ in the middle of Figure \ref{figHnbh} accounts for 4 percent, while the pattern $N^4$ on the right is extremely rare, with a frequency of 0.02 percent. There were also two patterns with a single neighbor, which arise when the blue piece is removed from $N^1$ and $N^2.$ They have frequency of 7 and 2 percent, respectively.  The average number of neighbors is 2.34.  Thus although the attractor is connected and contains many cycles, its degree of connectivity is rather low.

A few edges of the neighborhood graph can be found by inspection of Figure \ref{figHnbh}. Subpiece 2 on the right of the grey standard piece of $N^2$ has the same neighborhood $N^2$ as the grey piece. Thus there is a loop labelled 2 from $N^2$ to itself. Subpiece 1 on the left of the grey piece of $N^3$ and of $N^4$ has neighborhood $N^3,$ indicating a loop with label 1 at $N^3,$ and an edge with label 1 from $N^4$ to $N^3.$  The upper subpiece 3 of $N^3$ and $N^4$ has neighborhood $N^1.$ So there are edges with label 3 from $N^3$ and $N^4$ to $N^1.$ Altogether we have 90 edges for 30 neighborhoods and 3 mappings. As we have seen, edges describe the replacement of a piece by a subpiece -- that is, zooming in at the subpiece given by the edge label.

The loop at vertex $N^3$ is a reason for the fairly large frequency of this neighborhood. For an explanation of the extremely small frequency of $N^4,$ we have to see the graph. It turns out that $N_v=N^4$ only if the last letters of $v$ are $123(231)^k$ with $k\ge 1.$

\section{Applications of the neighborhood graph}  \label{out}
{\bf The interior as core of self-similarity.}  We started by asking how self-similar sets can be really self-similar when the global attractor $A$ and its magnifications appear so different.  For attractors with open set condition, we defined the interior $I(A).$ In Theorem \ref{intepie}c we proved that within each piece of $I(A)$ there are similar copies of all other pieces of the interior.
Thus $I(A)$ is truly self-similar. The magnifications inside $I(A)$ reveal the geometry of the attractor. For finite type attractors, the neighborhood graph is a scheme to produce these magnifications and determine their frequency.

The Sierpi\'nski triangle is the simplest example, with its three vertices as boundary sets.  When we magnify at any place outside the vertices, we get always the same photos, although in different order. If we magnify at a vertex, we get the very special views with the vertex. So we have to exclude vertices to get to the core of self-similarity. In other examples, the boundary to be excluded can be very large.

In order to discuss consequences, we identify the neighborhood of a piece with a local view of that piece in a certain window. To this end we would have to extend the neighbor concept, from neighbors which intersect the reference piece to neighbors which intersect the window, as indicated in \cite[Theorem 2.3]{EFG3}. The details do not fit into the present note. We confine ourselves to a brief non-technical review of related concepts in the literature and of work in progress.\vspace{1ex}

{\bf Frequencies of neighborhoods as discrete tangent distribution.}
The question which views are typical in a self-similar set with OSC was answered with the introduction of the tangential measure distribution \cite{Ba01,Gr95}. Even for much more general fractal measures, it was shown that continuous magnification around almost every point $x$ will lead to a collection of limit figures which does not depend anymore on $x$ \cite{Ga11,H10}.  Furstenberg \cite{Fu08} suggested the following terminology for a fractal $A$ and a window set $W$ in $\RR^d.$  The window is the unit cube or unit ball. It represents our screen or photo paper. A miniset or photo of $A$ has the form $g(A)\cap W$ where $g$ is a similitude. Usually $g$ is expanding with large expansion factor, so that only a small piece of $A$ will be mapped into our window $W.$  A microset of $A$ is a limit of minisets where the expansion factor tends to infinity. 

For a self-similar set $A$ with OSC and finite type, every microset is a miniset. That is, the tangential views arising in the limit of infinitely large magnification are already obtained on finite scale. Now Theorem \ref{intepie}c shows that even without finite type, only with OSC, every miniset of $I(A)$ occurs as a copy in arbitrary small pieces of $I(A),$ and thus is also a microset.  \emph{For the interior $I(A)$ of a self-similar set with OSC and finite type, the collections of minisets and of microsets coincide.}  The collection of minisets of $I(A)$ is invariant under magnification, and agrees with its own tangent. It is a linear object, like a line in elementary geometry. Linear objects deserve mathematical attention.

When we identify minisets with neighborhoods (which has to be worked out), there is only a finite number of different minisets, up to similarity.  We have shown that $I(A)$ is a countable union of pieces, but there are only finitely many pieces up to similarity. Thus a finite collection of photos will give all information about $A.$\vspace{1ex}

{\bf The neighborhood graph describes the magnification flow.}
Hochman \cite{H10} explained how the distribution of tangent measures develops from the successive magnification of neighborhoods of points in $A.$ For finite type attractors, and piece-centered minisets, our neighborhood graph discretizes this magnification flow. Following a path of edges in the neighborhood graph, we replace pieces by subpieces and subsubpieces.  Since we are in a finite graph, we can go on forever. 

The point is that we magnify ad infinitum with finitely many data. This can be realized practically, with one sufficiently accurate global picture of $A.$ Every local view is a finite union of copies of $A,$ described by certain neighbor maps. We just add together the shifted and rotated copies of the basic picture.  Only the list of neighbor maps and the neighbor graph is required. This is a virtual magnification. There are no problems with numerical accuracy.

We use discrete steps, scaling with $r^k,\ k\in\ZZ .$ Interpolation of continuous factors is possible, but not considered here.
When we zoom in, we have to select the subpiece which we want to see. A standard random method is to select each piece with probability $1/m.$  This defines a Markov chain between neighborhoods, with transition matrix  $\frac1m S,$ where $S$ is the substitution matrix in \eqref{subst}. The eigenvector $p=(p_1,...,p_K)$ of the leading eigenvalue $\sigma=m$ of $S$ is the stationary distribution of the Markov chain. It is the left eigenvector of the eigenvalue 1 of $\frac1m S.$ The right eigenvector is constant since we have a stochastic matrix. This justifies our definition of pattern frequencies.  

An interesting fact is that we can also zoom out indefinitely, just walking backwards the edge paths in the neighborhood graph. If our piece comes near the boundary of the set $I(A),$ we use Theorem \ref{intepie}c to find a tiny piece with similar neighborhood in the middle of $I(A),$ in order to find superpieces. In practice, we do not need the tiny piece.  Zooming out is virtually, like zooming in.  In each step we consult our graph to decide the neighborhood type of the superpiece of the current piece. 

Sometimes there is no choice.  The rare neighborhood $N^4$ in Figure \ref{figHnbh} is always piece 1 in its superpiece, which must be piece 3 in its supersuperpiece.  It depends on the number of edges directed towards the current neighborhood $N,$ as a vertex of the neighborhood graph.  If there are edges from several vertices $N^1,...,N^k$ to $N,$ a choice has to be made. The standard random method is to take the edge from $N^j$ with probability  $q_j=p_j/\sum_{n=1}^k p_n .$ Here $p_j$ denotes the stationary probability of the neighborhood $N^j$ in the Markov chain.

The Markov chain model is closed-in-itself, like a compact manifold. We can zoom out to see larger and larger scenery, only to arrive at some picture which we obtained by zooming in to very small detail.  The scenery is not periodic, however, because of the possible branchings at every step.  In the theory of aperiodic tilings and patterns, `branched manifolds' and `spaces of tilings' with a complicated solenoidal topology have been used as models \cite{AP}.  Our graph models are not manifolds or topological spaces.  They are easier to handle and provide a combinatorial way to look at the subject.\vspace{1ex}

{\bf Description of the examples by new parameters.} Besides visualization, the method opens new ways of describing the attractors. This is demonstrated here for the examples of Figures \ref{figA} and \ref{figB}.  Instead of Figure \ref{figC}, we consider the modification $v_1=1, v_2=1+i,$ and $v_3=-1$ of Figure \ref{figA}. The attractor will be called $A'.$ Beside the number of neighbors, of edges in the neighbor graph, and of neighborhoods, Table \ref{tab3} shows the minimum and maximum number of neighbors in a neighborhood. Average numbers of neighbors and frequencies are determined with the stationary distribution of the Markov chain of neighborhoods.  

\begin{table}[h] \label{tab3}       
\begin{tabular}{p{6mm}p{11mm}p{7mm}p{6mm}p{7.5mm}p{3.5mm}p{12.5mm}p{5.5mm}p{8mm}p{4mm}p{14mm}p{23mm}}
\hline
Name&neighbor&graph&&\# of&nbs.&in Nbh.&Freq.&of \# &nbs.& (\%) &Leading Nbh.\\
&vertices&edges&Nbh.&min&max&average&1&2&3&$>\!\frac23$ max&Frequencies \\
\hline
A & 23 & 41 & 88 & 1&6&2.98&11&21&38&11&12, 10, 7 \\
A' & 21 & 45 & 333 & 1&13&2.97&9&31&38& 0.2&13, 12, 8 \\
B & 82& 171 & 6291 & 1&14&2.56&10&35&47&0.0002&19, 10, 10\\
\hline
\end{tabular}
\caption{Some new parameters describing the structure of the fractal attractors $A,B$ from Section \ref{exa} and the modification $A'$ of $A.$ nbs.=neighbors, Nbh.=neighborhood, Freq.=frequency (percentage).}
\end{table}

The average number of neighbors is between 2.5 and 3 in all examples. The maximum number of neighbors in a neighborhood is 6, 13, and 14, respectively, although the total number of neighbors is much larger. Ten percent of the neighborhoods have a single neighbor, 21 up to 35 percent have two neighbors and about 40 percent have three neighbors, accounting together for 70 to 90 percent of the patterns.  The frequency of neighborhoods with more than $\frac23$ of the maximal number of neighbors is quite different: 11 percent for $A,$ 0.2 percent for $A1,$ and only 0.0002 percent for $B$ although these were more than thousand patterns.

In all examples, a few patterns have high frequency. The three largest frequencies, given in Table \ref{tab3}, sum up to 30-40 \% .
Approximately half of the neighborhoods belong to one of the five or six leading types.  The structure of these most frequent patterns is simple, with at most 3 neighbors. 

The well-connected pieces with many neighbors are most interesting and most influential for the global structure.  There can be many of these patterns, but they are rare, as the total frequency 0.0002 \% for attractor $B$ indicates.  The neighborhood graph makes it possible to look up these interesting parts of the attractors which would be otherwise difficult to find.
\vspace{1ex}

{\bf Outlook.}
This is exploratory research. While measure and dimension theory of fractals have been developed to a high degree of completeness \cite{BSS,Bar,BP,Fal,Fra}, the study of the topological structure is a wide open field.  Even simple examples can be hard to imagine. The present stage of computing technology opens new chances.  To this end we need a computational fractal geometry. Here we introduced the neighbor graph and the neighborhood graph as tools.

We studied only very simple examples in the plane with expansion factor 2, integer translations, and rotations by multiples of $90^o.$ The tools become more helpful for complicated IFS.  The number of neighborhoods can become very large, so that construction of the whole graph has to be replaced by sampling of neighborhoods.  It is necessary to pass from self-similar sets \eqref{hut} to graph-directed systems, and from plane to space. 

However, it is also possible to completely abandon similitudes and affine maps, and take the two neighbor automata as starting point \cite{EFG4}.  We can define the topology of $A$ by a neighbor automaton, and produce metric realizations by physical principles, for instance minimization of energy. This requires new ideas and bears the chance of approaching fractals in nature. \vspace{1ex}

{\bf Acknowledgement.} The author thanks Manuel Moran for correcting his first version of Theorem 3.1.


\begin{thebibliography}{10}

\bibitem{Tet24}
K.~Allabergenova, M.~Samuel, and A.~Tetenov.
\newblock Intersections of the pieces of self-similar dendrites in the plane.
\newblock {\em Chaos, Solitons and Fractals}, 182:114805, 2024.

\bibitem{AP}
J.E. Anderson and I.F. Putnam.
\newblock Topological invariants for substitution tilings and their associated
  ${C}^*$-algebras.
\newblock {\em Ergodic Theory and Dynam. Systems}, 18:509--537, 1998.

\bibitem{Baake2013}
M.~Baake and U.~Grimm.
\newblock {\em Aperiodic Order, Vol. 1: A mathematical invitation}.
\newblock Cambridge University Press, Cambridge, 2013.

\bibitem{Ba01}
C.~Bandt.
\newblock Local geometry of fractals given by tangent measure distributions.
\newblock {\em Monatshefte Math.}, 133:265--280, 2001.

\bibitem{EFG3}
C.~Bandt.
\newblock Elementary fractal geometry. 3. {C}omplex {P}isot factors imply
  finite type.
\newblock {\em Discrete Comput. Geom.}, 2024.

\bibitem{EFG4}
C.~Bandt.
\newblock Elementary fractal geometry. 4. {A}utomata-generated topological
  spaces.
\newblock {\em Communications in Math.}, 33(2):Paper no. 4, 2025.

\bibitem{EFG5}
C.~Bandt and M.F. Barnsley.
\newblock Elementary fractal geometry. 5. {W}eak separation is strong
  separation.
\newblock {\em arXiv:2404.04892}, 2024.

\bibitem{BHR}
C.~Bandt, N.V. Hung, and H.~Rao.
\newblock On the open set condition for self-similar fractals.
\newblock {\em Proc. Amer. Math. Soc.}, 134:1369--1374, 2005.

\bibitem{EFG1}
C.~Bandt and D.~Mekhontsev.
\newblock Elementary fractal geometry. {N}ew relatives of the {S}ierpi{\'n}ski
  gasket.
\newblock {\em Chaos: An Interdisciplinary Journal of Nonlinear Science},
  28(6):063104, 2018.

\bibitem{BM09}
C.~Bandt and M.~Mesing.
\newblock Self-affine fractals of finite type.
\newblock In {\em Convex and fractal geometry}, volume~84 of {\em Banach Center
  Publ.}, pages 131--148. Polish Acad. Sci. Inst. Math., Warsaw, 2009.

\bibitem{BSS}
B.~B{\'a}r{\'a}ny, K.~Simon, and B.~Solomyak.
\newblock {\em Self-similar and self-affine sets and measures}, volume 276 of
  {\em Mathematical Surveys and Monographs}.
\newblock American Mathematical Society, 2023.

\bibitem{Bar}
M.~F. Barnsley.
\newblock {\em Fractals everywhere}.
\newblock Academic Press, 2nd edition, 1993.

\bibitem{BP}
C.J. Bishop and Y.~Peres.
\newblock {\em Fractal sets in probability and analysis}.
\newblock Cambridge University Press, Cambridge, 2017.

\bibitem{DKV}
P.~Duvall, J.~Keesling, and A.~Vince.
\newblock The {Hausdorff} dimension of the boundary of a self-similar tile.
\newblock {\em J. London Math. Soc.}, 61:748--760, 2000.

\bibitem{Fal}
K.~J. Falconer.
\newblock {\em Fractal geometry: mathematical foundations and applications}.
\newblock J. Wiley \& sons, 3 edition, 2014.

\bibitem{Feng16}
D.-J. Feng.
\newblock The topology of polynomials with bounded integer coefficients.
\newblock {\em J. Eur. Math. Soc.}, 18:181--193, 2016.

\bibitem{Fra}
J.M. Fraser.
\newblock {\em {A}ssouad dimension and fractal geometry}.
\newblock Cambridge University Press, 2020.
\newblock arXiv 2005.03763.

\bibitem{Fu08}
H.~Furstenberg.
\newblock Ergodic fractal measures and dimension conservation.
\newblock {\em Ergodic Theory Dynam. Systems}, 28(2):205--222, 2008.

\bibitem{Ga11}
M.~Gavish.
\newblock Measures with uniform scaling scenery.
\newblock {\em Ergodic Theory Dynam. Systems}, 31(1):33--48, 2011.

\bibitem{Gi86}
W.J. Gilbert.
\newblock The fractal dimension of sets derived from complex bases.
\newblock {\em Canad. Math. Bull.}, 29(4):495--500, 1986.

\bibitem{Gr95}
S.~Graf.
\newblock On {B}andt's tangential distribution for self-similar measures.
\newblock {\em Monatsh. Math.}, 120:223--246, 1995.

\bibitem{GS}
B.~Gr\"{u}nbaum and G.C. Shephard.
\newblock {\em Patterns and Tilings}.
\newblock Freeman, New York, 1987.

\bibitem{HR22}
K.~E. Hare and A.~Rutar.
\newblock Local dimensions of self-similar measures satisfying the finite
  neighbor condition.
\newblock {\em Nonlinearity}, 35:4876--4904, 2022.

\bibitem{H10}
M.~Hochman.
\newblock Dynamics on fractals and fractal distributions.
\newblock {\em arXiv:1008.3731}, 2010.

\bibitem{Lal97}
S.P. Lalley.
\newblock $\beta$-expansions with deleted digits for {P}isot numbers $\beta$.
\newblock {\em Trans. Amer. Math. Soc.}, 349(11):4355--4365, 1997.

\bibitem{MW}
R.D. Mauldin and S.C. Williams.
\newblock Hausdorff dimension in graph-directed constructions.
\newblock {\em Trans. Amer. Math. Soc.}, 309:811--829, 1988.

\bibitem{M}
D.~Mekhontsev.
\newblock {IFS} tile finder, version 2.60.
\newblock \url{https://ifstile.com}, 2021.

\bibitem{Mo1}
M.~Moran.
\newblock Dynamical boundary of a self-similar set.
\newblock {\em Fund. Math.}, 160:1--14, 1999.

\bibitem{Mo}
P.A.P. Moran.
\newblock Additive functions of intervals and hausdorff measure.
\newblock {\em Math. Proc. Cambridge Phil. Soc.}, 42:15--23, 1946.

\bibitem{NW}
S.-M. Ngai and Y.~Wang.
\newblock Hausdorff dimension of self-similar sets with overlaps.
\newblock {\em J. London Math. Soc.}, 63:655--672, 2001.

\bibitem{Q}
M.~Queff{\' e}lec.
\newblock {\em Substitution dynamical systems - spectral analysis}, volume 1294
  of {\em Lecture Notes in Mathematics}.
\newblock Springer, Berlin, 1987.

\bibitem{rutar23}
A.~Rutar.
\newblock Geometrical and combinatorial properties of self-similar multifractal
  measures.
\newblock {\em Ergodic Theory Dyn. Syst.}, 43:2028--2072, 2023.
\newblock arXiv:2008.00197v3.

\bibitem{Sch}
A.~Schief.
\newblock Separation properties for self-similar sets.
\newblock {\em Proc. Amer. Math. Soc.}, 122:111--115, 1994.

\bibitem{SW}
R.S. Strichartz and Y.~Wang.
\newblock Geometry of self-affine tiles 1.
\newblock {\em Indiana Univ. Math. J.}, 48:1--24, 1999.

\bibitem{T89}
W.P. Thurston.
\newblock Groups, tilings, and finite state automata.
\newblock {A}MS Colloquium Lectures, Boulder, CO, 1989.

\bibitem{TZ20}
J.~Thuswaldner and S.~Zhang.
\newblock On self-affine tiles whose boundary is a sphere.
\newblock {\em Trans. Amer. Math. Soc.}, 373(1):491--527, 2020.

\end{thebibliography}
\end{document}